\documentclass[aos,MSNbibl,nameyear,dvips]{arximspdf}
\usepackage{dcolumn}
\usepackage{graphicx}

\doi{10.1214/11-AOS885}
\volume{39}
\issue{4}
\pubyear{2011}
\firstpage{1827}
\lastpage{1851}

\makeatletter

\newcolumntype{d}[1]{D{.}{.}{#1}}

\newtheorem{theorem}{Theorem}
\newtheorem{lemma}{Lemma}
\newproclaim{Remark}{Remark}

\newcommand{\mpl}{{\mathrm{MPL}}}
\newcommand{\pr}{\operatorname{pr}}

\newcommand{\Normal}{\operatorname{Normal}}

\newcommand{\T}{{\mathrm{T}}}
\newcommand{\sgn}{\operatorname{sgn}}
\newcommand{\diag}{\operatorname{diag}}
\newcommand{\SE}{\mathrm{SE}}
\newcommand{\wh}{\widehat}

\newcommand{\logit}{\operatorname{logit}}

\newcommand{\boldalpha}{\bolds{\alpha}}
\newcommand{\boldbeta}{\bolds{\beta}}
\newcommand{\boldgamma}{\bolds{\gamma}}
\newcommand{\boldnu}{\bolds{\nu}}
\newcommand{\boldomega}{\bolds{\omega}}

\newcommand{\boldtheta}{\bolds{\theta}}
\newcommand{\tbbeta}{\bolds{\beta}}
\newcommand{\tbnu}{\bolds{\nu}}

\makeatother

\begin{document}
\begin{frontmatter}

\title{Estimation and variable selection for generalized additive partial linear models}
\runtitle{Generalized additive partial linear models}

\begin{aug}
\author[A]{\fnms{Li} \snm{Wang}\thanksref{t1}\ead[label=e1]{lilywang@uga.edu}},
\author[B]{\fnms{Xiang} \snm{Liu}\thanksref{t2}\ead[label=e2]{xliu@bst.rochester.edu}},
\author[B]{\fnms{Hua} \snm{Liang}\corref{}\thanksref{t3}\ead[label=e3]{hliang@bst.rochester.edu}} and
\author[C]{\fnms{Raymond J.} \snm{Carroll}\thanksref{t4}\ead[label=e4]{carroll@stat.tamu.edu}}
\runauthor{Wang, Liu, Liang and Carroll}
\affiliation{University of Georgia, University of Rochester,
University of Rochester and
Texas A\&M University}
\address[A]{L. Wang\\
Department of Statistics\\
University of Georgia\\
Athens, Georgia 30602\\
USA\\
\printead{e1}}
\address[B]{X. Liu\\
H. Liang\\
Department of Biostatistics\\
\quad and Computational Biology\\
University of Rochester\\
Rochester, New York 14642\\
USA\\
\printead{e2}\\
\hphantom{E-mail: }\printead*{e3}}
\address[C]{R. J. Carroll\\
Department of Statistics\\
Texas A\&M University\\
College Station, Texas 77843-3143\\
USA\\
\printead{e4}}
\end{aug}

\thankstext{t1}{Supported by NSF Grant DMS-09-05730.}
\thankstext{t2}{Supported by a Merck Quantitative Sciences Fellowship Program.}
\thankstext{t3}{Supported by NSF Grants DMS-08-06097 and DMS-10-07167.}
\thankstext{t4}{Supported by a grant from the National Cancer Institute
(CA57030) and by Award Number KUS-CI-016-04, made by King Abdullah
University of Science and Technology (KAUST).}

\received{\smonth{3} \syear{2010}}
\revised{\smonth{2} \syear{2011}}

%
\begin{abstract}
We study generalized additive partial linear models, proposing the use
of polynomial spline smoothing for estimation of nonparametric
functions, and deriving quasi-likelihood based estimators for the
linear parameters. We establish asymptotic normality for the estimators
of the parametric components. The procedure avoids solving large
systems of equations as in kernel-based procedures and thus results in
gains in computational simplicity. We further develop a~class of
variable selection procedures for the linear parameters by employing a
nonconcave penalized quasi-likelihood, which is shown to have an
asymptotic oracle property. Monte Carlo simulations and an empirical
example are presented for illustration.
\end{abstract}

%

\end{frontmatter}

\section{Introduction}\label{SECintroduction}

Generalized linear models (GLM), introduced by \citet{nelwed72}
and systematically summarized by \citet{mccnel89}, are a powerful
tool to analyze the relationship
between a~discrete response variable and covariates. Given a link
function, the GLM expresses the relationship between the dependent
and independent variables through a linear functional form. However,
the GLM and associated methods may not be flexible enough when
analyzing complicated data generated from biological and
biomedical research. The generalized additive model (GAM), a
generalization of the GLM that replaces linear components by
a~sum of smooth unknown functions of predictor variables, has been
proposed as an alternative and has been used widely
[\citet{hastib90}, \citet{woobook06}]. The generalized
additive partially
linear model (GAPLM) is a realistic, parsimonious candidate when
one believes that the relationship between the dependent variable
and some of the covariates has a parametric form, while
the relationship between the dependent variable and the remaining
covariates may not be linear. GAPLM enjoys the simplicity
of the GLM and the flexibility of the GAM because it
combines both parametric and nonparametric components.

There are two possible approaches for estimating the parametric
component and the nonparametric components in a GAPLM. The first
is a combination of kernel-based backfitting and local scoring,
proposed by \citet{bujhastib89} and detailed by
\citet{hastib90}. This method may need to solve a
large system of equations [\citet{yuparkmammen08}], and also
makes it difficult to introduce a penalized function for
variable selection as given in Section \ref{SECvariableselection}. The
second is an application of the marginal
integration approach [\citet{linnie95}] to the
nonparametric component of the generalized partial linear models.
They treated the summand of additive terms as a
nonparametric component, which is then estimated as a multivariate
nonparametric function. This strategy may still suffer from the
``curse of dimensionality'' when the number of additive terms is
not small [\citet{hardleheutmammensperlich04}].

The kernel-based backfitting and marginal integration approaches
are computationally expensive.
\citet{MarxEilers98}, \citet{rupwancar03}
and \citet{WoodJASA04} studied penalized regression splines,
which share most of the practical benefits of smoothing spline
methods, combined with ease of use and reduction of the
computational cost of backfitting GAMs. Widely used R/Splus
packages \textsf{gam} and \textsf{mgcv}
provide a convenient implementation in practice. However, no
theoretical justifications are available for these procedures in the
additive case.
See \citet{lirup08} for recent work in the one-dimensional case.

In this paper, we will use polynomial splines to estimate the
nonparametric components. Besides asymptotic theory, we develop a
flexible and convenient estimation procedure for GAPLM. The use of
polynomial spline smoothing in generalized nonparametric
models goes back to \citet{stone86}, who first obtained the rate of
convergence of the polynomial spline estimates for the generalized
additive model. \citet{stone94} and \citet{huang98} investigated
polynomial spline estimation for the generalized functional ANOVA
model. More recently, \citet{xueyangSS06} studied estimation of
the additive coefficient model for a continuous response variable
using polynomial spline methods. Our models
emphasize possibly non-Gaussian responses, and combine both
parametric and nonparametric components through a link function.
Estimation is achieved through maximizing the
quasi-likelihood with polynomial
spline smoothing for the nonparametric functions. The convergence
results of the maximum likelihood estimates for the nonparametric
parts in this article are similar to those for regression
established by \citet{xueyangSS06}. However, it is very
challenging to establish asymptotic normality in our general
context, since it cannot be viewed simply as an orthogonal
projection, due to its nonlinear structure. To the best of our
knowledge, this is the \textit{first attempt} to establish
asymptotic normality of the estimators for the parametric
components in GAPLM. Moreover, polynomial spline smoothing is a
global smoothing method, which approximates the unknown functions
via polynomial splines characterized by a linear combination of
spline basis. After the spline basis is chosen, the coefficients
can be estimated by an efficient one-step procedure of maximizing
the quasi-likelihood function. In
contrast, kernel-based methods, such as those reviewed above, in
which the maximization must be conducted repeatedly at every
data point or a grid of values, are more time-consuming. Thus, the
application of polynomial spline smoothing in the current context
is particularly computationally efficient compared to some of its
counterparts.

In practice, a large number of variables may be collected and some
of the insignificant ones should be excluded before forming a
final model. It is an important issue to select significant
variables for both parametric and nonparametric regression models;
see \citet{fanli06} for a comprehensive overview of variable
selection. Traditional variable selection procedures such as
stepwise deletion and subset selection may be extended to the
GAPLM. However, these are also computationally expensive because,
for each submodel, we encounter the challenges mentioned above.

To select significant variables in semiparametric models, \citet
{liliang08} adopted Fan and Li's (\citeyear{fanli01}) variable selection
procedures for parame\-tric models via nonconcave penalized
quasi-likelihood, but their models do not cover the GAPLM. Of course,
before developing justifiable
variable selec\-tion for the GAPLM, it is important to establish
asymptotic properties for the parametric components. In this
article, we propose a class of variable selection procedures for
the parametric component of the GAPLM and study the asymptotic
properties of the resulting estimator. We demonstrate how the rate
of convergence of the resulting estimate depends on the
regulariza\-tion parameters, and further show that the penalized
\mbox{quasi-likelihood}~esti\-mators perform asymptotically as an oracle procedure
for selecting the~model.

The rest of the article is organized as follows. In Section
\ref{SECmodels}, we introduce the GAPLM model. In Section
\ref{SECestimation}, we propose polynomial spline estimators via
a quasi-likelihood approach, and study the asymptotic properties
of the proposed estimators. In Section \ref{SECvariableselection}, we
describe the variable selection procedu\-res for the
parametric component, and then prove their statistical properties.
Simulation studies and an empirical example are presented in
Section~\ref{secexam}. Regularity conditions and the proofs of the
main results are presented in the \hyperref[app]{Appendix}.

\section{The models}\label{SECmodels}

Let $Y$ be the response variable, $\mathbf{X}=(
X_{1},\ldots,X_{d_1})^{\T}\in R^{d_1}$ and $\mathbf{Z}=(
Z_{1},\ldots,Z_{d_2})^{\T}\in R^{d_2}$ be the covariates. We
assume the conditional density of $Y$ given $(
\mathbf{X},\mathbf{Z}) =( \mathbf{x},\mathbf{z})
$ belongs to the exponential family
%
\begin{equation}\label{DEFloglikelihood}
f_{ Y| \mathbf{X},\mathbf{Z}}(  y|
\mathbf{x},%
\mathbf{z}) =\exp[ y\xi( \mathbf{x},\mathbf
{z}) -%
\mathcal{B}\{ \xi( \mathbf{x},\mathbf{z})
\} +%
\mathcal{C}( y) ]
\end{equation}
for known functions $\mathcal{B}$ and $\mathcal{C}$, where $\xi$
is the so-called natural parameter in parametric generalized
linear models (GLM), is related to the unknown mean response
by
\[
\mu( \mathbf{x},\mathbf{z}) =E( Y|\mathbf
{X}=\mathbf{x},%
\mathbf{Z}=\mathbf{z}) =\mathcal{B}^{\prime}\{ \xi
( \mathbf{x},\mathbf{z}) \} .
\]
In parametric GLM, the mean function $\mu$ is defined via a known
link function~$g$ by $g\{ \mu( \mathbf{x},\mathbf{z})
\}
=\mathbf{x}^{\T}\boldalpha+\mathbf{z^{\T}\boldbeta}$, where
$\boldalpha$ and $\boldbeta$ are parametric vectors to be
estimated. In this article, $g( \mu) $ is modeled as
an additive partial linear function
%
\begin{equation}\label{modelGAPLM}
g\{ \mu( \mathbf{x},\mathbf{z}) \}
=\sum_{k=1}^{d_1}\eta_{k}(
x_{k})+\mathbf{z}^{\T}\boldbeta,
\end{equation}
where $\boldbeta$ is a $d_2$-dimensional regression parameter,
$\{\eta_{k}\}_{k=1}^{d_1}$ are unknown and smooth
functions and $E\{\eta_{k
}( X_{k})\}=0$ for $1\leq k \leq d_1$ for
identifiability.

If the conditional variance function
$\operatorname{var}(Y|\mathbf{X}=\mathbf{x},\mathbf{Z}=\mathbf{z}) =\sigma^{2}V\{ \mu( \mathbf{x},\mathbf{z}) \}$
for some known positive function $V$, then estimation of the mean
can be achieved by replacing the conditional loglikelihood function
$\log\{f_{Y| \mathbf{X},\mathbf{Z}}( y|\allowbreak \mathbf{x},\mathbf{z})\}$ in
(\ref{DEFloglikelihood}) by a quasi-likelihood function $Q(
m,y) $, which satisfies
\[
\frac{\partial}{\partial m}Q( m,y) =\frac{y-
m}{V( m)}.
\]
The first goal of this article is to provide a simple method of
estimating $\boldbeta$ and $\{\eta
_{k}\}_{k=1}^{d_1}$ in model (\ref{modelGAPLM}) based on a
quasi-likelihood procedure [\citet{sevsta94}] with
polynomial splines. The second goal is to discuss how to select
significant parametric variables in this semiparametric framework.

\section{Estimation method} \label{SECestimation}

\subsection{Maximum quasi-likelihood} \label{subseclikelihood}
Let $( Y_{i},\mathbf{X}_{i},\mathbf{Z}_{i}) $,
$i=1,\ldots,n$, be independent copies of $(
Y,\mathbf{X},\mathbf{Z}) $. To avoid confusion, let $\eta
_{0}=\sum_{k=1}^{d_1}\eta_{0k}( x_{k})$ and
$\boldbeta_{0}$ be the true additive function and the true
parameter values, respectively.
{For simplicity, we assume that the covariate
$X_{k }$ is distributed on a compact interval $[a_{k },b_{k }]$},
$k =1,\ldots,d_1$, and without loss of generality, we take all
intervals $[a_{k },b_{k }]=[0,1]$, $k =1,\ldots,d_1$. Under
smoothness assumptions, the~$\eta_{0k}$'s can be well
approximated by spline functions. Let~$\mathcal{S}_{n}$ be the
space of polynomial splines on $[0,1]$ of order $r \geq1$. We
introduce a knot sequence with $J$ interior knots
\[
\xi_{-r+1
}=\cdots=\xi_{-1}=\xi_{0}=0<\xi_{1}<\cdots<\xi_{J}<1=\xi_{J+1}=\cdots=\xi_{J+r},
\]
where $J\equiv J_{n}$ increases when sample size $n$ increases,
where the precise order is given in condition (C5) in Section \ref{subsecassumption-result}.
According to \citet{stone85},~$\mathcal{S}_{n}$ consists of functions $\hbar$
satisfying:\vspace*{3pt}
\begin{longlist}
\item $\hbar$ is a polynomial of
{degree $r-1$} on each of the
subintervals $I_{j}=[ \xi_{j},\xi_{j+1}) $,
$j=0,\ldots,J_{n}-1$, $I_{J_n}=[ \xi_{J_{n}},1] $;
\item for $r \geq2$, $\hbar$ is $r-2$ times continuously
differentiable on $[0,1]$.\vspace*{3pt}
\end{longlist}
Equally-spaced knots are used in this article for simplicity of
proof. However other regular knot sequences can also be used, with
similar asymptotic results.

We will consider additive spline estimates $\widehat{\eta}$ of
$\eta_{0}$. Let $\mathcal{G}_{n}$ be the
collection of functions $\eta$ with the additive form $\eta(
\mathbf{x%
}) =\sum_{k=1}^{d_1}\eta_{k}(
x_{k})$, where each component function $\eta_{k }\in\mathcal
{S}%
_{n}$ and $\sum_{i=1}^{n}\eta_{k }( X_{ik }) =0$. We
seek a function $\eta\in\mathcal{G}_{n}$ and a value of $
\boldbeta$ that maximize the quasi-likelihood function
%
\begin{equation}\label{DEFquasilikelihood}
L( \eta,\boldbeta) =n^{-1}\sum_{i=1}^{n}Q[
g^{-1}\{ \eta( \mathbf{X}_{i})
+\mathbf{Z}_{i}^{\T}\boldbeta\} ,Y_{i}].
\end{equation}
For\vspace*{1pt} the $k$th covariate $x_{k}$, let $b_{j,k}( x_{k })
$ be the B-spline basis functions of order $r$. For any $\eta\in
\mathcal{G}_{n}$, write\vspace*{1pt} $\eta( \mathbf{x}) =\boldgamma
^{\T}\mathbf{b}( \mathbf{x})$, where $\mathbf{b}(
\mathbf{x}) =\{ b_{j,k }( x_{k }),
j=1,\ldots,J_{n}+r, k =1,\ldots,d_1\}^{\T}$
{is the collection of} the spline basis functions, and
$\boldgamma=\{ \gamma_{j,k }, j=1 ,\ldots,J_{n}+r, k =1,\ldots,d_1
\} ^{\T}$ is the spline coefficient vector. Thus, the
maximization problem in (\ref{DEFquasilikelihood}) is equivalent
to finding $\boldbeta$ and~$\boldgamma$ to maximize
%
\begin{equation}\label{DEFquasilikelihood1}
\ell(\boldgamma,\boldbeta)=n^{-1}\sum_{i=1}^{n}Q[
g^{-1}\{\boldgamma^{\T}\mathbf{b}( \mathbf{X}_{i}
) +%
\mathbf{Z}_{i}^{\T}\boldbeta\} ,Y_{i}].
\end{equation}
We denote\vspace*{2pt} the maximizer as $\widehat{\boldbeta}$ and $
\widehat{\boldgamma}=\{ \widehat{\gamma}_{j,k }, j=1
,\ldots,J_{n}+r, k =1,\ldots,d_1\} ^{\T}$. Then the spline
estimator of $\eta_{0}$ is $\widehat{\eta}(
\mathbf{x})=\widehat{\boldgamma}^{\T}\mathbf{b}(
\mathbf{x})$, and the centered spline component function
estimators are
\[
\widehat{\eta}_{k }( x_{k })
=\sum_{j=1}^{J_{n}+r}\widehat{\gamma}_{{j,k}}b_{j,k}
( x_{k }) -n^{-1}%
\sum_{i=1}^{n}\sum_{j=1 }^{J_{n}+r}\widehat{\gamma
}_{{j,k}}b_{j,k}( X_{ik }),\qquad k=1,\ldots,d_1.
\]
The above estimation approach can be easily implemented because
this approximation results in a generalized linear model. However,
theoretical justification for this estimation approach is very
challenging [\citet{huang98}].

Let $N_n=J_{n}+r-1$. We adopt the normalized B-spline space
$\mathcal{S}_{n}^{0}$ introduced in \citet{xueyangSS06} with the
following normalized basis
%
\begin{eqnarray}\label{DEFBjk}
B_{j,k}( x_{k}) =\sqrt{N_{n}}\biggl\{ b_{j+1,k}(
x_{k}) -\frac{E( b_{j+1,k})}{E(
b_{1,k})} b_{1,k}( x_{k}) \biggr\},\nonumber\\[-8pt]\\[-8pt]
&&\eqntext{1\leq j\leq
N_n,  {1}\leq k\leq d_{1},}
\end{eqnarray}
which is convenient for asymptotic analysis. Let $\mathbf{B}(
\mathbf{x}) =\{ B_{j,k }( x_{k }), j=1,\ldots,\allowbreak N_{n},\ k =1,\ldots,d_1\}
^{\T}$ and $\mathbf{B}_{i}= \mathbf{B}%
( \mathbf{X}_{i})$. Finding $(\boldgamma,
\boldbeta) $ that maximizes (\ref{DEFquasilikelihood1}) is
mathematically equivalent to finding $(\boldgamma,
\boldbeta)$ which maximizes
\[
n^{-1}\sum_{i=1}^{n}Q[
g^{-1}\{\mathbf{B}_{i}^{\T}\boldgamma
+\mathbf{Z}_{i}^{\T}\boldbeta\} ,Y_{i}].
\]
Then the spline estimator of $\eta_{0}$ is $\widehat{\eta}(
\mathbf{x})= \widehat{\boldgamma}^{\T}\mathbf{B}(
\mathbf{x}) $, and the centered spline estimators of the
component functions are
\[
\widehat{\eta}_{k }( x_{k }) =\sum_{{j=2}
}^{N_{n}}\widehat{\gamma}_{j,k}B_{j,k}( x_{k }) -n^{-1}%
\sum_{i=1}^{n}\sum_{j=2}^{N_{n}}\widehat{\gamma}_{j,k}B_{j,k}(
X_{ik }),\qquad k =1,\ldots,d_1.
\]
We show next that estimators of both the parametric and nonparametric
components have nice asymptotic properties.

\vspace*{-2pt}\subsection{Assumptions and asymptotic results}\label{subsecassumption-result}

Let $v$ be a positive integer and $\alpha\in( 0,1] $
such that
$p=v+\alpha>2$. Let $\mathcal{H}{(p)}$ be the collection of
functions $g$ on $[0,1]$ whose $v$th derivative,
$g^{(v) }$, exists and satisfies a Lipschitz
condition of order $\alpha$, $| g^{( v ) }(
m^{*}) -g^{( v )
}( m) | \leq C|
{m^{*}}-m| ^{\alpha}$,
for $0\leq{m^{*}},m\leq1$,
where $C$ is a positive constant. Following the notation of
\citet{carfangij97}, let $\rho_{\ell}( m)=
\{
{dg^{-1}( m)/dm}
\} ^{\ell}/ V\{ g^{-1}( m)
\} $ and
$q_{\ell}( m,y) =\partial^{\ell}/\partial
m^{\ell}Q\{ g^{-1}( m) ,y\} $, so that
\begin{eqnarray*}
q_{1}( m,y) &=&\partial/\partial m Q\{
g^{-1}( m) ,y\} =\{ y-g^{-1}( m)
\} \rho_{1}( m),\\
q_{2}( m,y) &=&\partial^{2}/\partial m^{2}Q\{
g^{-1}( m) ,y\} =\{ y-g^{-1}( m)
\} \rho_{1}^{\prime}( m) -\rho_{2}(
m).
\end{eqnarray*}

For simplicity of notation, write $\mathbf{T}=( \mathbf
{X},\mathbf{Z}
) $ and $\mathbf{A}^{\otimes2}=\mathbf{A} \mathbf{A}^{\T}$
for any matrix or vector $\mathbf{A}$. We make the following
assumptions:

\begin{longlist}[(C5)]
\item[(C1)] The function $\eta_{0}^{{\prime\prime}}( \cdot
)$ is continuous and each component function $\eta
_{0k}(\cdot)\in\mathcal{H}{(p)}$, $k =1,\ldots,d_1$.

\item[(C2)] The function $q_{2}( m,y) <0$ and
$c_{q}<| q_{2}^{\nu}( m,y) | <C_{q}$
($\nu=0,1$) for $m\in R$ and $y$ in the range of the response
variable.

\item[(C3)] The distribution of $\mathbf{X}$ is absolutely
continuous and its density $f$ is bounded away from zero and
infinity on $[0,1]^{d_1}$.

\item[(C4)] The random vector $\mathbf{Z}$ satisfies that for any
unit vector $ \boldomega\in R^{d_2}$
\[
c\leq\boldomega^{\T}E(
\mathbf{Z}^{\otimes2}|\mathbf{X}=\mathbf{x})\boldomega\leq
C.
\]

\item[(C5)] The number of knots $n^{1/( 2p) }\ll
N_{n}\ll n^{1/4}$.
\end{longlist}
\begin{Remark}\label{Remark1}
The smoothness condition in (C1) describes
a requirement on the best rate of convergence that the functions
$\eta_{0k} (\cdot)$'s can be approximated by functions in the
spline spaces. Condition (C2) is imposed to ensure the uniqueness
of the solution; see, for example, Condition 1a of \citet
{carfangij97} and Condition (i) of \citet{liliang08}. Condition (C3)
requires a boundedness condition on the covariates, which is
often assumed in asymptotic analysis of nonparametric regression
problems; see Condition 1 of \citet{stone85}, Assumption (B3)(ii)
of \citet{huang99} and Assumption (C1) of \citet{xueyangSS06}.
The boundedness assumption on the support can be replaced by a finite
third moment assumption, but this will add much extra complexity to the
proofs. Condition (C4) implies that the eigenvalues of
$E(\mathbf{Z}^{\otimes2}|\mathbf{X}=\mathbf{x})$ are
bounded away from $0$ and $\infty$. Condition~(C5) gives the rate
of growth of the dimension of the spline spaces relative to the
sample size.
\end{Remark}

For measurable functions $\varphi_{1}$, $\varphi_{2}$ on
$[ 0,1]^{d_1}$, define the empirical inner product and
the corresponding norm as
\[
\langle\varphi_{1},\varphi_{2}\rangle
_{n}=n^{-1}\sum_{i=1}^{n}\{
\varphi_{1}( \mathbf{X}_{i}) \varphi_{2}( \mathbf
{X}%
_{i}) \},\qquad \| \varphi\|
_{n}^{2}=n^{-1}\sum_{i=1}^{n}\varphi^{2}(
\mathbf{X}_{i}) .
\]
If $\varphi_{1}$ and $\varphi_{2}$ are $L^{2}$-integrable,
define the theoretical inner product and corresponding norm as
\[
\langle\varphi_{1},\varphi_{2}\rangle=E\{ \varphi
_{1}( \mathbf{X}) \varphi_{2}( \mathbf{X})
\},\qquad \| \varphi\| _{2}^{2}=E\varphi
^{2}( \mathbf{X}).
\]
Let $\| \varphi\| _{nk }^{2}$ and $\| \varphi
\|
_{2k }^{2}$ be the empirical and theoretical norm of $\varphi$
on $[ 0,1] $, defined by
\[
\| \varphi\| _{nk }^{2}=n^{-1}\sum_{i=1}^{n}\varphi
^{2}( X_{ik }),\qquad \| \varphi\|
_{2k }^{2}=E\varphi^{2}( X_{k }) =\int_{0}^{1}\varphi
^{2}( x_{k }) f_{k }( x_{k }) \,dx_{k },
\]
where $f_k(\cdot)$ is the density function of $X_k$.

Theorem \ref{THMconvergence} describes the rates of convergence of the
nonparametric parts.

\begin{theorem}
\label{THMconvergence} Under conditions \textup{(C1)--(C5)},
for $k =1,\ldots,d_1$, $\| \widehat{\eta
}-\eta_{0}\| _{2} =O_{P}\{ N_{n}^{1/2-p}+ (N_n/n)^{1/2}\}$;
$\| \widehat{\eta}-\eta_{0}\| _{n} =O_{P}\{
N_{n}^{1/2-p}+ (N_n/n)^{1/2}\}$;\break $\| \widehat{\eta}_{k
}-\eta_{0k }\| _{2k } =O_{P}\{ N_{n}^{1/2-p}+
(N_n/n)^{1/2}\}$ and\vspace*{1pt} $\| \widehat{\eta}_{k }-\eta_{0k
}\| _{nk } =O_{P}\{ N_{n}^{1/2-p}+ (N_n/ n)^{1/2}\}$.
\end{theorem}

Let $m_{0}( \mathbf{T}) =\eta_{0}( \mathbf{X}%
) +\mathbf{Z}^{\T}\boldbeta_{0}$ and define
%
\begin{equation} \label{DEFGammax}
\Gamma( \mathbf{x}) =\frac{E[  \mathbf{Z}%
\rho_{2}\{ m_{0}( \mathbf{T}) )\} |
\mathbf{X}=%
\mathbf{x}] }{E[  \rho_{2}\{ m_{0}(
\mathbf{T}%
) \} | \mathbf{X}=\mathbf{x}]
},\qquad
\widetilde{\mathbf{Z}}=\mathbf{Z}-\Gamma^{\mathrm{add}}(
\mathbf{X}),
\end{equation}
where
%
\begin{equation}\label{DEFaddproj}
\Gamma^{\mathrm{add}}(\mathbf{x})=\sum_{k=1}^{d_{1}}\Gamma_{k}(x_{k})
\end{equation}
is the\vspace*{1pt} projection of $\Gamma$ onto the Hilbert space of
theoretically centered additive functions with a norm
$\|f\|_{\rho_{2},m_0}^2=E[f(\mathbf{X})^2\rho_{2}\{m_0(\mathbf{T})\}]$.
To obtain asymptotic normality of the estimators in
the linear part, we further impose the conditions:
\begin{longlist}[(C8)]
\item[(C6)] The additive components in (\ref{DEFaddproj})
satisfy that $\Gamma_{k}(\cdot)\,{\in}\,\mathcal{H}{(p)},\,k\,{=}\,1,\ldots,d_1$.
\item[(C7)] For $\rho_{\ell}$, we have
\begin{eqnarray}
| \rho_{\ell}( m_{0}) | \leq C_{\rho
}\quad\mbox{and}\quad| \rho_{\ell}( m) -\rho
_{\ell}( m_{0}) | \leq C_{\rho}^{*}|
m-m_{0}| \nonumber\\
&&\eqntext{\mbox{for all }| m-m_{0}| \leq C_{m},
\ell=1,2.}
\end{eqnarray}

\item[(C8)] There exists a positive constant $C_{0}$, such that
$E[  \{ Y-g^{-1}( m_{0}(\mathbf{T}) ) \} ^{2}|\allowbreak \mathbf{T}]
\leq C_{0}$, almost surely.
\end{longlist}

The next theorem shows that the maximum quasi-likelihood estimator
of~$\boldbeta_{0}$ is root-$n$ consistent and asymptotically normal,
although the convergence rate of the nonparametric component $
\bolds{\eta}_{0}$ is of course slower than root-$n$.
\begin{theorem}
\label{THMnormality}
$\!\!\!\!\!$Under conditions \textup{(C1)--(C8)},
$ \sqrt{n}(\widehat{\boldbeta}\,{-}\,\boldbeta_{0})\,{\rightarrow}\,\Normal( 0,\bolds{\Omega}^{-1})$,
where $\bolds{\Omega}=E[ \rho_{2}\{ m_{0}(\mathbf{T}) \}\widetilde{\mathbf{Z}}^{\otimes2}]$.
\end{theorem}

The proofs of these theorems are given in the \hyperref[app]{Appendix}.

It is worthwhile pointing out that taking the additive structure of the
nuisance parameter into account leads to a smaller asymptotic variance
than that of the estimators which ignore the additivity
[\citet{YuLee10}]. \citet{CarrollMaityMammenYu09} had the same
observation for a special case with repeated measurement data when $g$
is the identity function.

\section{Selection of significant parametric variables}\label{SECvariableselection}

In this section, we develop variable selection procedures for the
parametric component of the GAPLM. We study the asymptotic
properties of the resulting estimator, illustrate how the rate of
convergence of the resulting estimate depends on the
regularization parameters, and further establish the oracle
properties of the resulting estimate.

\subsection{Penalized likelihood}

Building upon the quasi-likelihood given in
(\ref{DEFquasilikelihood}), we define the penalized
quasi-likelihood as
%
\begin{equation}\label{DEFpenalized-quasilikelihood1}
\mathcal{L}(\eta,\boldbeta)=\sum_{i=1}^{n}Q[ g^{-1}\{
\eta( \mathbf{X}_{i})
+\mathbf{Z}_{i}^{\T}\boldbeta\}
,Y_{i}]-n\sum_{j=1}^{d_2}p_{\lambda_{j}}(
| \beta_{j}| ) ,
\end{equation}
where\vspace*{1pt} $p_{\lambda_{j}}( \cdot) $ is a prespecified
penalty function with a regularization parameter~$\lambda_{j}$.
The penalty functions and regularization parameters in
(\ref{DEFpenalized-quasilikelihood1}) are not necessarily the
same for all $j$. For example, we may wish to keep scientifically
important variables in the final model, and therefore do not want
to penalize their coefficients. In practice, $\lambda_{j}$ can be
chosen by a data-driven criterion, such as cross-validation (CV)
or generalized cross-validation [GCV, \citet{cravenwahba79}].

Various penalty functions have been used in variable selection for
linear regression models, for instance, the $L_{0}$ penalty, in
which $p_{\lambda_{j}}( | \beta| )
=0.5\lambda_{j}^{2}I( |\beta|\ne0) $. The
traditional best-subset variable selection can be viewed as a
penalized least squares with the $L_{0}$ penalty because
$\sum_{j=1}^{d_2}I( | \beta_{j}| \ne
0)$ is essentially the number of nonzero regression
coefficients in the model. Of course, this procedure has two well
known and severe problems. First, when the number of covariates is
large, it is computationally infeasible to do subset selection.
Second, best subset variable selection suffers from high
variability and instability [\citet{breiman96}, \citet{fanli01}].

The Lasso is a regularization technique for simultaneous
estimation and variable selection [\citet{tib96}, \citet{zou06}] that
avoids the drawbacks of the best subset selection. It can be
viewed as a penalized least squares estimator with the $L_{1}$
penalty, defined by $p_{\lambda_{j}}( | \beta|
) =\lambda_{j}|\beta|$. \citet{frafri93}
considered bridge regression with an $L_{q}$ penalty, in which
$p_{\lambda_{j}}( | \beta| ) =\lambda
_{j}|\beta|^{q}$ ($0<q<1$). The issue of selection of the
penalty function has been studied in depth by a variety of
authors. For example, \citet{fanli01} suggested using the SCAD
penalty, defined by
\begin{eqnarray}
p_{\lambda_{j}}^{\prime}(\beta)=\lambda_{j}\biggl\{ I(\beta\le
\lambda_{j})+\frac{(a\lambda_{j}-\beta)_{+}}{(a-1)\lambda
_{j}}I(\beta>\lambda_{j})\biggr\}\nonumber\\
&&\eqntext{\mbox{for some }a>2
\mbox{ and } \beta>0,}
\end{eqnarray}
where $p_{\lambda_{j}}( 0) =0$, and $\lambda_{j}$ and
$a$ are two tuning parameters. \citet{fanli01} suggested using
$a=3.7$, which will be used in Section \ref{secexam}.

Substituting $\eta$ by its estimate in
(\ref{DEFpenalized-quasilikelihood1}), we obtain a penalized
likelihood
%
\begin{equation} \label{EQLpbeta}
\mathcal{L}_{P}(\boldbeta)=\sum_{i=1}^nQ[g^{-1}\{\mathbf
{B}_{i}^{\T}
\widehat{\boldgamma}+\mathbf{Z}_{i}^{\T}\boldbeta\}
,Y_{i}]-n\sum_{j=1}^{d_2}p_{\lambda_{j}}(
| \beta_{j}| ).
\end{equation}
Maximizing $\mathcal{L}_{P}(\boldbeta)$ in (\ref{EQLpbeta})~yields a maximum penalized likelihood
estimator~$\widehat{\boldbeta}^{\mpl}$. The theorems established below demonstrate
that~$\widehat{\boldbeta}^{\mpl}$ performs asymptotically as well as an oracle estimator.

\subsection{Sampling properties}

We next show that with a proper choice of $\lambda_{j}$, the
maximum penalized likelihood estimator
$\widehat{\boldbeta}^{\mpl}$ has an asymptotic oracle property. Let
$\boldbeta_{0}=(\beta_{10},\ldots,\beta
_{d_20})^{\T}=(\boldbeta_{10}^{\T},\boldbeta_{20}^{\T})^{\T}$,
where $\boldbeta_{10}$ is assumed to consist of all nonzero
components\vspace*{1pt} of $\boldbeta_{0}$ and $\boldbeta_{20}=\mathbf{0}$
without loss of generality. Similarly we write
$\mathbf{Z}=(\mathbf{Z}_{1}^{\T},\mathbf{Z}_{2}^{\T})^{\T}$. Denote
$w_{n}=\max_{1\le j\le{d_{2}}}\{\vert p_{\lambda_{j}}^{\prime\prime} ( \vert\beta_{j0}\vert) \vert,\allowbreak\beta_{j0}\ne0\} $
and
%
\begin{equation}\label{DEFan}
a_{n}=\max_{1\le j\le d_2}\{ \vert
p_{\lambda_{j}}^{\prime}(\vert\beta_{j0}\vert) \vert,\beta
_{j0}\ne0\}.
\end{equation}
\begin{theorem}\label{THMrootn}
Under the regularity conditions given in
Section \ref{subsecassumption-result}, and if $a_{n}\to0$ and
$w_{n}\to0$ as $n\to\infty$,
then there exists a local maximizer $\widehat{\boldbeta}^{\mpl}$
of $\mathcal{L}_{P}(\boldbeta)$ defined in (\ref{EQLpbeta}) such
that its rate of convergence is $O_{P}(n^{-1/2}+a_{n})$, where
$a_{n}$ is given in (\ref{DEFan}).
\end{theorem}

Next, define $\bolds{\xi}_{n}=\{p_{\lambda_{1}}^{\prime}(|\beta
_{10}|) \sgn(\beta_{10}),\ldots,p_{\lambda_{s}}^{\prime}(|\beta
_{s0}|) \sgn(\beta_{s0})\}^{\T}$ and a diagonal matrix
$\bolds{\Sigma} _{\lambda}=\diag\{p_{\lambda_{1}}^{\prime
\prime}(|\beta_{10}|),\ldots,p_{\lambda_{s}}^{\prime\prime
}(|\beta_{s0}|)\}$, where $s$ is the number of nonzero components
of $\boldbeta_{0}$. Define $\mathbf{T}_{1}=(
\mathbf{X},\mathbf{Z}_{1}) $ and $m_{0}(
\mathbf{T}_{1}) =\eta_{0}( \mathbf{X})
+\mathbf{Z}_{1}^{\T} \boldbeta_{10}$, and further let
\[
\Gamma_{1}( \mathbf{x}) =\frac{E[  \mathbf
{Z}%
_{1}\rho_{2}\{ m_{0}( \mathbf{T}_{1}) \}
| \mathbf{X}=\mathbf{x}] }{E[  \rho
_{2}\{ m_{0}( \mathbf{T}_{1} ) \} |
\mathbf{X}=\mathbf{x}]},\qquad \widetilde{\mathbf{Z}}_{1}
=\mathbf{Z}_{1}-\Gamma_{1}^{\mathrm{add}}( \mathbf{X}),
\]
where $\Gamma_{1}^{\mathrm{add}}$ is the projection of $\Gamma_{1}$ onto
the Hilbert space of theoretically centered additive functions
with the norm $\|f\|_{\rho_{2},m_0}^2$.
\begin{theorem} \label{THMoracle} Suppose that\vspace*{2pt} the
regularity conditions given in Section~\ref{subsecassumption-result}
hold, and that $\liminf_{n\to\infty}\liminf_{\beta_{j}\to0^{+}}\lambda_{jn}^{-1}
p_{\lambda_{jn}}^{\prime}(|\beta_{j}|)>0$. If $\sqrt{n}\lambda_{jn}\to\infty$ as $n\to\infty$, then the
root-$n$ consistent estimator $\widehat{\boldbeta}^{\mpl}$ in Theorem \ref{THMrootn} satisfies
$\widehat{\boldbeta}_{2}^{\mpl}\!=\mathbf{0}$, and $\sqrt{n}(\bolds{\Omega}_{s}+
\bolds{\Sigma}_{\lambda})\{\widehat{\boldbeta}_{1}^{\mpl}\!-\boldbeta_{10}+(\bolds{\Omega}_{s}+
\bolds{\Sigma}_{\lambda})^{-1}\bolds{\xi}_{n}\}\,{\rightarrow}\,\Normal(\mathbf{0},\bolds{\Omega}_s)$,
where $\bolds{\Omega}_{s}=[\rho_{2}\{m_{0}(\mathbf{T}_{1})\}\widetilde{\mathbf{Z}}_{1}^{\otimes2}]$.
\end{theorem}

\subsection{Implementation}

As pointed out by \citet{liliang08}, many penalty functions,
including the $L_{1}$ penalty and the SCAD penalty, are irregular
at the origin and may not have a second derivative at some points.
Thus, it is often difficult to implement the Newton--Raphson
algorithm directly. As in \citet{fanli01}, \citet{hunli05},
we approximate the penalty function locally by a quadratic function at
every step in the iteration such that the
Newton--Raphson algorithm can be modified for finding the solution
of the penalized likelihood. Specifically, given an initial
value~$\bolds{\beta}^{(0)}$ that is close to the maximizer of the
penalized likelihood function, the penalty $p_{\lambda_{j}}(
| \beta_{j}| ) $ can be locally approximated by
the quadratic function as $ \{p_{\lambda_{j}}( | \beta
_{j}| )\} ^{\prime}=p_{\lambda_{j}}^{\prime}(
| \beta_{j}| ) \sgn(\beta_j) \approx
\{p_{\lambda_{j}}^{\prime}( \vert\beta_{j}^{(0)}\vert) / \vert
\beta_{j}^{(0)}\vert\} \beta_{j}, $ when $\beta_{j}^{(0)}$ is
not very close to 0; otherwise, set $\widehat{\beta}_{j}=0$. In
other words, for $\beta_{j}\approx\beta_{j}^{(0)}$, $
p_{\lambda_{j}}( | \beta_{j}| ) \approx
p_{\lambda_{j}}(\vert\beta_{j}^{(0) }\vert) + (1/2) \{p_{\lambda
_{j}}^{\prime}( \vert\beta_{j}^{(0)}\vert) / \vert\beta
_{j}^{(0)}\vert\} ( \beta_{j}^{2}-\beta_{j}^{(0)2})$. For
instance, this local quadratic approximation for the $L_{1}$
penalty yields
\[
| \beta_{j}| \approx(1/2)\bigl\vert\beta_{j}^{(0)}\bigr\vert
+
(1/2)\beta_{j}^{2}/ \bigl\vert\beta_{j}^{(0)}\bigr\vert %
\qquad\mbox{for } \beta_{j}\approx\beta_{j}^{(0)}.
\]

\vspace*{6pt}

\textit{Standard error formula for}
$\widehat{\boldbeta}^{\mpl}$. We follow\vspace*{1pt} the approach in \citet
{liliang08} to derive a sandwich formula for the estimator
$\widehat{\boldbeta}^{\mpl}$. Let
\begin{eqnarray*}\ell^{\prime
}(\boldbeta)&=&\frac{\partial\ell
(\widehat{\boldgamma},\boldbeta)}{\partial\boldbeta},\qquad  {\ell
}^{\prime\prime}(\boldbeta)=\frac{\partial^{2} \ell
(\widehat{\boldgamma},\boldbeta)}{\partial\boldbeta\,
\partial\boldbeta^{\T}}; \\
\bolds{\Sigma} _{\lambda}(\boldbeta)&=&\diag\biggl\{
\frac{p_{\lambda_{1}}^{\prime}(|\beta_{1}|)}{|\beta
_{1}|},\ldots,\frac{p_{\lambda
_{{d_{2}}}}^{\prime}(|\beta
_{{d_{2}}}|)} {|\beta
_{{d_{2}}}|}\biggr\}.
\end{eqnarray*}
A sandwich formula is given by
\begin{eqnarray*}
\widehat{\operatorname{cov}}( \widehat{\boldbeta}^{\mpl})&=&\{n
{\ell}%
^{\prime\prime}( \widehat{\boldbeta}^{\mpl}) -n\bolds{\Sigma}
_{\lambda
}( \widehat{\boldbeta}^{\mpl}) \} ^{-1}\widehat{\operatorname{cov}}%
\{ {\ell}^{\prime}( \widehat{\boldbeta}^{\mpl}) \} \\
&&{} \times\{n {\ell}^{\prime\prime} (
\widehat{\boldbeta}^{\mpl}) -n\bolds{\Sigma}_{\lambda}(
\widehat{\boldbeta}^{\mpl}) \} ^{-1}.
\end{eqnarray*}
Following conventional techniques that arise in the likelihood
setting, the above sandwich formula can be shown to be a consistent
estimator and will be shown in our simulation study to have good
accuracy for moderate sample sizes.\vspace*{8pt}

\textit{Choice of $\lambda_{j}$\textup{'}s}. The unknown
parameters $(\lambda_{j})$ can be selected using data-driven
approaches, for example, generalized cross validation as proposed
in \citet{fanli01}. Replacing $\boldbeta$ in
(\ref{DEFquasilikelihood1}) with its estimate
$\widehat{\boldbeta}^{\mpl}$, we maximize $\ell
(\boldgamma,\widehat{\boldbeta}^{\mpl})$ with respect to
$\boldgamma$. The solution\vadjust{\goodbreak} is denoted by
$\widehat{\boldgamma}^{\mpl}$, and the corresponding estimator of
$\eta_{0}$ is defined as
%
\begin{equation}\label{DEFetaMPL}
\widehat{\bolds{\eta}}^{\mpl}(\mathbf{x})=
(\widehat{\boldgamma}
^{\mpl})^{\T}\mathbf{B}(\mathbf{x}).
\end{equation}
Here the GCV statistic is defined by
\[
\operatorname{GCV}(\lambda_{1},\ldots,\lambda
_{d_2})=\frac{\sum_{i=1}^{n}D[
Y_{i},g^{-1}\{ \widehat{\bolds{\eta}}^{\mpl}(
\mathbf{X}_{i})
+\mathbf{Z}_{i}^{\T}\widehat{\boldbeta}^{\mpl}\} ] }{
n\{1-e(\lambda_{1},\ldots,\lambda
_{{d_{2}}})/n\}^{2}},
\]
where $ e(\lambda_{1},\ldots,\lambda
_{{d_{2}}})=\operatorname{tr}[ \{
\ell^{\prime\prime}(\widehat{\boldbeta}^{\mpl})
-n\bolds{\Sigma} _{\lambda}(\widehat{\boldbeta}^{\mpl}) \}
^{-1}\ell^{\prime\prime}(\widehat{\boldbeta}^{\mpl}) ] $ is the
effective number of parameters and $D(Y,\mu)$ is the deviance of
$Y$ corresponding to fitting with $\bolds{\lambda}$. The
minimization problem over a $d_{2}$-dimensional space is
difficult. However, \citet{liliang08} conjectured that the
magnitude of $\lambda_{j}$ should be proportional to the standard
error of the unpenalized maximum pseudo-partial likelihood
estimator of $ \beta_{j}$. Thus, we suggest taking $\lambda
_{j}=\lambda\SE( \widehat{\beta}_{j})$ in practice, where
$\SE(\widehat{\beta}_{j}) $ is the estimated standard error
of~$\widehat{\beta}_{j}$, the unpenalized likelihood estimate
defined in Section \ref{SECestimation}. Then the minimization
problem can be reduced to a one-dimensional problem, and the
tuning parameter can be estimated by a grid search.

\vspace*{3pt}\section{Numerical studies}\vspace*{3pt}
\label{secexam}

\subsection{A simulation study}\label{subsecsim}

We simulated $100$ data sets consisting of $n=100$, $200$ and
$400$ observations, respectively, from the GAPLM:
%
\begin{equation}\label{eq51}
\logit\{\pr(Y=1)\}=\eta_1(X_1)+ \eta_2(X_2) + \mathbf{Z}^{\T
}\boldbeta,
\end{equation}
where
\begin{eqnarray*}
\eta_1(x)&=&\sin(4 \pi x),\\
\eta_2(x)&=&10 \{\exp(-3.25x) +
4\exp(-6.5x)+3\exp(-9.75x)\}
\end{eqnarray*}
and the true parameters $\boldbeta
=(3,1.5,0,0,0,0,2,0)^{\T}$.
$X_1$ and $X_2$ are independently uniformly distributed on $[0,1]$.
$Z_1$ and $Z_2$ are normally distributed with mean $0.5$ and variance
$0.09$. The random vector
$(Z_1,\ldots,Z_6,X_1,X_2)$ has an autoregressive structure with correlation
coefficient $\rho=0.5$.

In order to determine the number of knots in the approximation, we
performed a simulation with 1,000 runs for each sample size. In each
run, we fit, without any variable selection procedure, all possible
spline approximations with 0--7 internal knots for each nonparametric
component. The internal knots were equally spaced quantiles from the
simulated data. We recorded the combination of the numbers of knots
used by the best approximation, which had the smallest prediction
error~(PE), defined as
%
\begin{equation}\label{eqpredictionerror}
\mathrm{PE} =\frac1n \sum^{n}_{i=1}\bigl\{\logit^{-1}(\mathbf{B}_{i}^{\T}
\widehat{\boldgamma}+\mathbf{Z}_i^{\T}\widehat{\boldbeta})
-\logit^{-1}\bigl(\eta(\mathbf{X}_i)+\mathbf{Z}_i^{\T}{\boldbeta}\bigr)\bigr\}^2.
\end{equation}
$(2,2)$ and $(5,3)$ are most frequently chosen for sample sizes $100$
and $400$, respectively. These combinations were used in the
simulations for the variable selection procedures.\vadjust{\goodbreak}

The proposed selection procedures were applied to this model and
B-splines were used to approximate the two nonparametric functions. In
the simulation and also the empirical example in Section \ref{subsecreal}, the estimates from ordinary logistic regression were used
as the starting values in the fitting procedure.

To study model fit, we also defined model error (ME) for the parametric
part by
%
\begin{equation}\label{eqmodelerror}
\mathrm{ME}(\wh{\boldbeta})=
(\wh{\boldbeta}-\boldbeta)^{\T}E(ZZ^{\T})(\wh{\boldbeta
}-\boldbeta).
\end{equation}
The relative model error is defined as the ratio of the model error
between the fitted model using variable selection methods and using
ordinary logistic regression.

The simulation results are reported in Table \ref{ta1}, in which
the columns labeled with ``$C$'' give the average number of the five
zero coefficients correctly set to~$0$, the columns labeled with
``$I$'' give the average number of the three nonzero coefficients
incorrectly set to $0$, and the columns labeled with ``MRME'' give
the median of the relative model errors.

\begin{table}
\caption{Results from the simulation study in Section
\protect\ref{subsecsim}. $C$, $I$ and MRME stand for the average
number of the five
zero coefficients correctly set to $0$,
the average number of the three nonzero coefficients incorrectly set to $0$,
and the median of the relative model errors. The model errors are
defined in
(\protect\ref{eqmodelerror})}\label{ta1}
\begin{tabular*}{\tablewidth}{@{\extracolsep{\fill}}lcccc@{}}
\hline
$\bolds{n}$ & \multicolumn{1}{c}{\textbf{Method}}
& \multicolumn{1}{c}{$\bolds{C}$} & \multicolumn{1}{c}{$\bolds{I}$}
& \multicolumn{1}{c@{}}{\textbf{MRME}} \\
\hline
100& ORACLE &5\hphantom{.00}&0\hphantom{.00}&0.27\\
& SCAD&4.29&0.93&0.60\\
& Lasso&3.83&0.67&0.51\\
& BIC&4.53&0.95&0.54\\
[4pt]
400& ORACLE&5\hphantom{.00}&0\hphantom{.00}&0.33\\
& SCAD&4.81&0.27&0.49\\
& Lasso&3.89&0.10&0.67\\
& BIC&4.90&0.35&0.46\\
\hline
\end{tabular*}
\vspace*{-6pt}
\end{table}

Summarizing Table \ref{ta1}, we conclude that BIC performs the best in
terms of correctly identifying zero coefficients, followed by SCAD and
LASSO. On the other hand, BIC is also more likely to set nonzero
coefficients to zero, followed by SCAD and LASSO. This indicates that
BIC most aggressively reduce the model complexity, while LASSO tends to
include more variables in the models. SCAD is a useful compromise
between these two procedures. With an increase of sample sizes, both
SCAD and BIC nearly perform as if they had Oracle property. The MRME
values of the three procedures are comparable. Results of the cases not
depicted here have characteristics similar to those shown in Table \ref{ta1}.
Readers may refer to the online supplemental materials.

\begin{figure}

\includegraphics{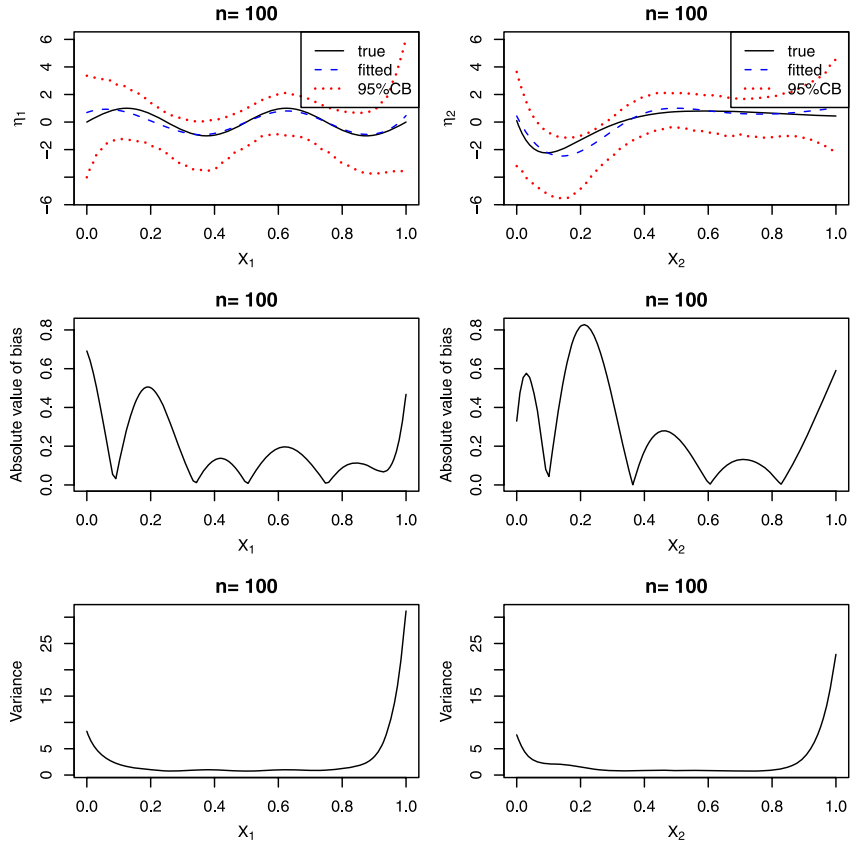}

\caption{The mean, absolute value of the bias and variance of the
fitted nonparametric functions when $n=100$ and $\rho=0.5$ [the left
panel for $\eta_1(x_1)$ and the right for $\eta_2(x_2)$]. $95\%$ CB
stands for the $95\%$ confidence band.} \label{fignonpar100}
\end{figure}

\begin{figure}

\includegraphics{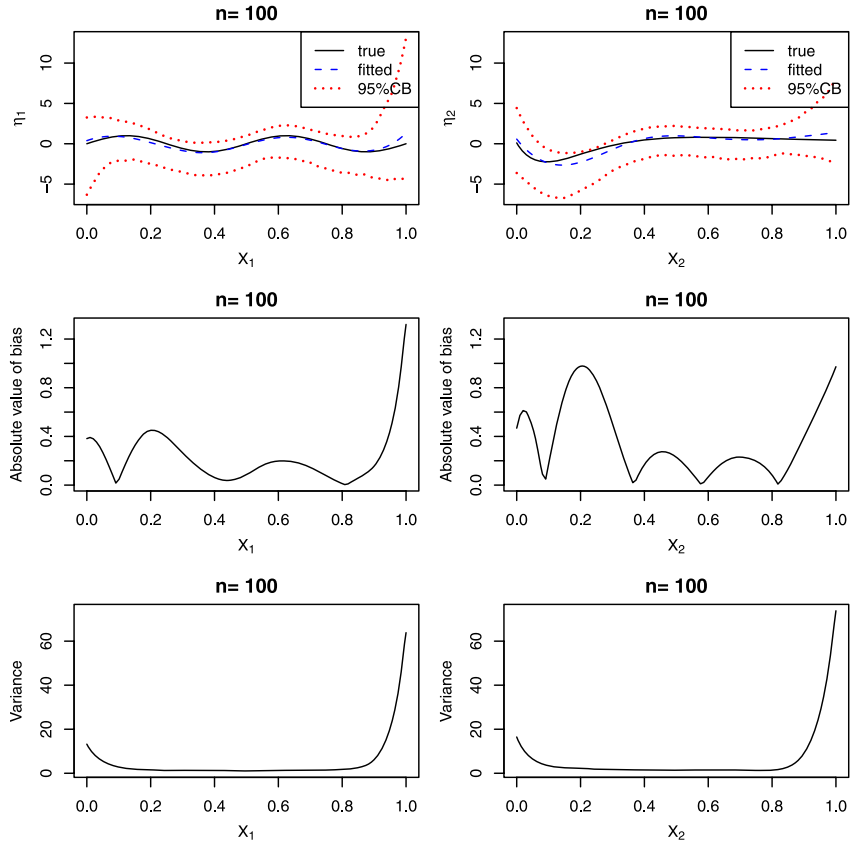}

\caption{The mean, absolute value of the bias and variance of the
fitted nonparametric functions when $n=100$ and $\rho=0.7$. The left
panel is for $\eta_1(x_1)$ and the right panel is for~$\eta_2(x_2)$.
Here $95\%$ CB stands for the $95\%$ confidence band.}
\label{fignonpar100v13}
\vspace*{-3pt}
\end{figure}

We also performed a simulation with correlated covariates.
We generated the response $Y$ from model (\ref{eq51}) again but with
$\boldbeta=(3.00, 1.50, 2.00)$. The covariates $Z_1$, $Z_2$, $X_1$ and
$X_2$ were marginally normal with mean zero and variance $0.09$. In
order, $(Z_1, Z_2, X_1, X_2)$ had autoregressive correlation
coefficient $\rho$, while $Z_3$ is Bernoulli with success probability $0.5$.
We considered two scenarios: (i) moderately correlated covariates
($\rho=0.5$) and
(ii) highly correlated ($\rho=0.7$) covariates.
We did 1,000 simulation runs
for each case with sample sizes $n=100, 200$ and $400$. From our
simulation, we
observe that the estimator becomes more unstable when the correlation among
covariates is higher. In scenario (i), all simulation runs
converged. However, there were $6$, $3$ and $7$ cases of
nonconvergence over
the 1,000 simulation runs for sample sizes $100, 200$ and $400$,
respectively,
in scenario (ii). In addition, the variance and bias of the fitted
functions in
scenario (ii) were much larger than those in scenario (i), especially
on the boundaries
of the covariates' support. This can be observed in Figures
\ref{fignonpar100} and \ref{fignonpar100v13}, which present the mean,
%
%
absolute value of bias and variance
of the fitted nonparametric
functions for
$\rho=0.5$ and $\rho=0.7$ with sample size $n=100$. Similar results are
obtained for sample sizes $n=200$ and $400$, but are not given here.

\subsection{An empirical example}\label{subsecreal}

We now apply the GAPLM and our variable selection procedure to a
data set from the Pima Indian diabetes study [\citet{smitheverhartetc88}].
This data set is obtained from the UCI Repository of
Machine Learning Databases, and is selected from a larger data set
held by the National Institutes of Diabetes and Digestive and
Kidney Diseases. All patients in this database are Pima Indian
women at least 21 years old and living near Phoenix, Arizona. The
response $Y$ is the indicator of a positive test for diabetes.
Independent variables from this data set include: $\mathit{NumPreg}$, the
number of pregnancies; $\mathit{DBP}$, diastolic blood pressure (mmHg);
$\mathit{DPF}$, diabetes pedigree function; $\mathit{PGC}$, the plasma glucose
concentration after two hours in an oral glucose tolerance test;
$\mathit{BMI}$, body mass index [weight in kg$/$(height in m)$^2$]; and $\mathit{AGE}$
(years). There are in total $724$~complete observations in this
data set.

In this example, we explore the impact of these covariates
on the probability of a positive test. We first fit the data set
using a linear logistic regression model: the estimated results
are listed in the left panel of Table \ref{tabpima}. These
results indicate that $\mathit{NumPreg}$, $\mathit{DPF}$, $\mathit{PGC}$ and $\mathit{BMI}$ are
statistically significant, while $\mathit{DBP}$ and $\mathit{AGE}$ are not
statistically significant.

\begin{table}
\tabcolsep=0pt
\caption{Results for the Pima study. Left panel: estimated values,
associated standard errors and $P$-values by using GLM. Right
panel: Estimates, associated standard errors using the GAPLM with
the proposed variable selection procedures} \label{tabpima}
\begin{tabular*}{\tablewidth}{@{\extracolsep{\fill}}ld{2.3}d{1.3}d{2.3}d{1.3}ccc@{}}
\hline
& \multicolumn{4}{c}{\textbf{GLM}} &\multicolumn{3}{c@{}}{\textbf{GAPLM}}
\\[-4pt]
& \multicolumn{4}{c}{\hspace*{-2pt}\hrulefill} &\multicolumn{3}{c@{}}{\hrulefill} \\
& \multicolumn{1}{c}{\textbf{Est.}} & \multicolumn{1}{c}{\textbf{s.e.}} & \multicolumn{1}{c}{$\bolds{z}$ \textbf{value}}
& \multicolumn{1}{c}{\textbf{Pr}$\bolds{(}\mbox{$\bolds{>}$}\bolds{|z|)}$} & \multicolumn{1}{c}{\textbf{SCAD (s.e.)}}
& \mbox{\textbf{LASSO (s.e.)}} & \multicolumn{1}{c@{}}{\textbf{BIC (s.e.)}} \\
\hline
NumPreg& 0.118& 0.033& 3.527& 0&0 (0)\hphantom{000000..}&\hphantom{$-$}0.021 (0.019)&0 (0)\hphantom{000000..}\\
DBP& -0.009& 0.009& -1.035& 0.301&0 (0)\hphantom{000000..}&$-$0.006 (0.005)&0 (0)\hphantom{000000..}\\
DPF& 0.961& 0.306& 3.135& 0.002&0.958 (0.312)&\hphantom{$-$}0.813 (0.262)&0.958 (0.312)\\
PGC& 0.035& 0.004& 9.763& 0&0.036 (0.004)&\hphantom{$-$}0.034 (0.003)&0.036 (0.004)\\
BMI& 0.091& 0.016& 5.777& 0 & & \\
AGE& 0.017& 0.01& 1.723& 0.085 & & \\
\hline
\end{tabular*}
\vspace*{-3pt}
\end{table}

\begin{figure}

\includegraphics{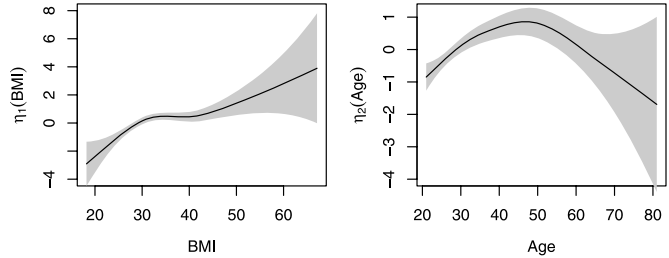}

\caption{The patterns of the nonparametric functions of BMI and Age
(solid lines) with${}\pm{}$s.e. (shaded areas) using the R function, gam,
for the Pima study.} \label{figpima}
\vspace{-3pt}
\end{figure}

However, a closer investigation shows that the effect of $\mathit{AGE}$ and
$\mathit{BMI}$ on the logit transformation of the probability of a positive
test may be nonlinear, see Figure \ref{figpima}. Thus, we employ
the following GAPLM for this data analysis,
%
\begin{eqnarray}\label{eqf1}
\logit\{P{(Y=1)}\}&=&\eta_0 + \beta_1 \mathit{NumPreg} + \beta
_2 \mathit{DBP} + \beta_3 \mathit{DPF} \nonumber\\[-8pt]\\[-8pt]
&&{} + \beta_4 \mathit{PGC} + \eta_1(\mathit{BMI}) + \eta_2(\mathit{AGE}).
\nonumber
\end{eqnarray}
Using B-splines to approximate $\eta_1(\mathit{BMI})$ and $\eta_2(\mathit{AGE})$,
we adopt $5$-fold cross-validation to select knots and find that the
approximation with no internal knots performs well for the both
nonparametric components.

We applied the proposed variable selection procedures to the model
(\ref{eqf1}), and the estimated coefficients and their
standard errors are listed in the right panel of
Table \ref{tabpima}.
Both SCAD and BIC suggest that $\mathit{DPF}$ and $\mathit{PGC}$ enter the model, whereas
$\mathit{NumPreg}$ and $\mathit{DBP}$ are suggested not to enter. However, the LASSO
suggests an inclusion of $\mathit{NumPreg}$ and $\mathit{DBP}$. This may be because LASSO
admits many variables in general,
as we observed in the simulation studies. The nonparametric estimators of
$\eta_1(\mathit{BMI})$ and $\eta_2(\mathit{AGE})$, which are obtained by using the
SCAD-based procedure, are similar to the solid lines\vadjust{\eject} in
Figure \ref{figpima}. It is worth pointing that the effect of
$\mathit{AGE}$ on the probability of a positive test shows a concave
pattern, and women whose age is around 50 have the highest
probability of developing diabetes. Importantly, the linear
logistic regression model does not reveal this significant effect.

It is interesting that the variable $\mathit{NumPreg}$ is statistically
insignificant when we fit the data using GAPLM with the proposed
variable selection procedure, but shows a statistically
significant impact when we use GLM. One might reasonably
conjecture that this phenomenon might be due to model
misspecification. To test this, we conducted a simulation as
follows. We generated the response variables using the estimates
and functions obtained by GAPLM with the SCAD. Then we fit a GLM
for the generated data set. We repeated the generation and fitting
procedures 5,000 times and found that $\mathit{NumPreg}$ is identified
positively significant $67.42\%$ percent of the time at level
$0.05$ in the GLMs. For $\mathit{DBP}$, $\mathit{DPF}$, $\mathit{PGC}$, $\mathit{BMI}$ and $\mathit{AGE}$, the
percentages that they are identified as statistically significant
at the level $0.05$ are $4.52\%$, $90.36\%$, $100\%$ and $99.98\%$
and $56.58\%$, respectively. This means that $\mathit{NumPreg}$ can incorrectly
enter the model, with more than $65\%$ probability, when a wrong model
is used, while $\mathit{DBP}$, $\mathit{DPF}$, $\mathit{PGC}$, $\mathit{BMI}$ and $\mathit{AGE}$ seem correctly to
be classified as insignificant and significant covariates even with
this wrong GLM model.

\section{Concluding remarks}

We have proposed an effective polynomial spline technique for the
GAPLM, then developed variable selection procedures to identify which
linear predictors should be included in the final model fitting.
The contributions we made to the existing literature can be
summarized in three ways: (i)~the procedures are computationally
efficient, theoretically reliable, and intuitively appealing; (ii) the
estimators of the linear components, which are often of primary
interest, are asymptotically normal; and (iii) the variable
selection procedure for the linear components has an asymptotic oracle
property. We believe that our approach can be extended to the
case of longitudinal data [\citet{LinCarroll06}], although the
technical details are by no means straightforward.

An important question in using GAPLM in practice is which covariates
should be included in the linear component. We suggest proceeding as
follows. The continuous covariates are put in the
nonparametric part and the discrete covariates in the parametric part.
If the estimation results show that some of the continuous covariate effects
can be described by certain parametric forms such as a linear form,
either by formal testing or by visualization, then a new model can be
fit with those continuous covariate effects moved to the parametric
part. The procedure can be iterated several times if needed. In this
way, one can take full advantage of the flexible exploratory analysis
provided by the proposed method. However, developing a more efficient
and automatic criterion warrants future study.
It is worth pointing\vadjust{\eject} out the proposed procedure may be instable for
high-dimensional data, and may encounter collinear problems. Addressing
these challenging questions is part of ongoing work.

\begin{appendix}\label{app}
\section*{Appendix}

Throughout the article, let $\| \cdot\| $ be the
Euclidean norm and $\| \varphi\| _{\infty
}=\sup_{m}| \varphi( m) | $ be the supremum
norm of a function $\varphi$ on $[
0,1] $. For any matrix~$\mathbf{A}$, denote its $L_{2}$ norm as
$%
\| \mathbf{A}\| _{2}=\sup_{\| \mathbf{x}\|
\neq0}%
{\|\mathbf{A}\mathbf{x}\| }/{\| \mathbf{x}
\|
}$, the largest eigenvalue.

\vspace*{-3pt}\subsection{Technical lemmas}

In the following, let $\mathcal{F}$
be a class of measurable functions. For probability\vspace*{1pt} measure $Q$, the
$L_{2}( Q)$-norm of a function $f\in\mathcal{F}$ is
defined by $(\int|f|^2 \,dQ)^{1/2}$.
According to \citet{vanderVaartWellner96}, the $\delta$-covering number
$\mathcal{N}( \delta, \mathcal{F},L_{2}( Q)
) $ is the smallest value of~$\mathcal{N}$
for which there exist functions $f_{1},\ldots,f_{\mathcal{N}}$, such
that for each $f\in\mathcal{F}$, $\| f-f_{j}\| \leq
\delta$ for some $j\in\{ 1,\ldots,\mathcal{N}\} $. The
$\delta$-covering\vspace*{1pt} number with bracketing $\mathcal{N}_{[
\cdot] }( \delta,\mathcal{F},L_{2}( Q) )
$ is the smallest value of $\mathcal{N}$ for which there exist
pairs\vspace*{-1pt} of functions $\{ [ f_{j}^{L},f_{j}^{U}]
\} _{j=1}^{\mathcal{N}}$ with $\|
f_{j}^{U}-f_{j}^{L}\| \leq\delta$, such that for each
$f\in\mathcal{F}$, there~is~a~$j\in\{1,\ldots,\mathcal{N}\}$ such that $f_{j}^{L}\leq f\leq
f_{j}^{U}$. The $\delta$-entropy with \mbox{bracketing}~is~defined as
$\log\mathcal{N}_{[\cdot]}( \delta,\mathcal{F},L_{2}(Q))$. Denote
$\mathcal{J}_{[\cdot]}( \delta,\mathcal{F},L_{2}( Q) ) =
\int_{0}^{\delta}\sqrt{1+\log\mathcal{N}_{[\cdot]}(\varepsilon,\mathcal{F},L_{2}( Q) )} \,d\varepsilon$.
Let $Q_{n}$ be the empirical measure of $Q$. Denote $G_{n}=
\sqrt{n}( Q_{n}-Q) $ and $\| G_{n}\|
_{\mathcal{F}} ={\sup_{f\in\mathcal{F}}}| G_{n}f| $
for any measurable class of functions $\mathcal{F}$.

We state several preliminary lemmas first, whose proofs are included in
the supplemental materials. Lemmas
\ref{LEMExpection-entropy}--\ref{LEMWn} will be used to prove
the remaining lemmas and the main results. Lemmas
\ref{LEMbetatilde} and \ref{LEMthetatilde-thetahat} are used to
prove \mbox{Theorems~\ref{THMconvergence}--\ref{THMrootn}}.
\begin{lemma}[{[Lemma 3.4.2 of \citet
{vanderVaartWellner96}]}]\label{LEMExpection-entropy}
Let $M_{0}$ be a finite positive constant. Let $\mathcal{F}$ be a
uniformly bounded class of measurable functions such that
$Qf^{2}<\delta^{2}$ and $\| f\| _{\infty}<M_{0}$.
Then
\[
E_{Q}^{*}\| G_{n}\| _{\mathcal{F}}\leq
C_{0}\mathcal{J}_{[\cdot]
}( \delta,\mathcal{F},L_{2}( Q) ) \biggl\{
1+\frac{%
\mathcal{J}_{[\cdot] }( \delta
,\mathcal{F},L_{2}( Q) ) }{\delta
^{2}\sqrt{n}}M_{0}\biggr\},
\]
where $C_{0}$ is a finite constant independent of $n$.
\end{lemma}
\begin{lemma}[{[Lemma A.2 of \citet{huang99}]}]\label{LEMEntropy}
For any $\delta>0$, let
\[
\Theta_{n}=\{ \eta( \mathbf{x}) +\mathbf{z}^{\T}
\bolds\boldbeta%
;\|\boldbeta-\boldbeta_{0}\| \leq\delta,\eta\in
\mathcal{G}_{n},\| \eta-\eta_{0}\| _{2}\leq\delta
\}.
\]
Then, for any $\varepsilon\leq\delta$, $ \log
\mathcal{N}_{[ \cdot] }( \delta,\Theta
_{n},L_{2}( P) ) \leq c N_n \log(\delta
/\varepsilon)$.
\end{lemma}

For simplicity, let
%
\begin{equation}\label{DEFDiWn}
\mathbf{D}_{i}=(\mathbf{B}_{i}^{\T},\mathbf{Z}_{i}^{\T}),\qquad
\mathbf{W}_{n}=n^{-1}\sum_{i=1}^{n}
\mathbf{D}_{i}^{\T}\mathbf{D}_{i}.
\end{equation}
\begin{lemma}\label{LEMWn}
$\!\!\!\!\!$Under conditions \textup{(C1)--(C5)}, for the above random
matrix~$\mathbf{W}_{n}$, there exists a positive constant $C$ such that
$\| \mathbf{W}_{n}^{-1}\| _{2}\leq C$, a.s.
\end{lemma}

According to a result of de Boor [(\citeyear{deBoor01}), page 149], for any
function $g\in\mathcal{H}{(p)}$ with $p<r-1$, there exists a
function $\widetilde{g}\in\mathcal{S}_{n}^{0}$, such that
$\| \widetilde{g}-g\| _{\infty}\leq CN_{n}^{-p}$,
where $C$ is some fixed positive constant. For $\eta_{0}$
satisfying (C1), we can find $\widetilde{\boldgamma}=\{
\widetilde{\gamma}_{j,k }, j=1 ,\ldots,N_{n}, k =1,\ldots,d_1\}
^{\T} $ and an\vspace*{1pt} additive spline function $\widetilde{\eta}=\widetilde
{%
\boldgamma}^{\T}\mathbf{B}( \mathbf{x}) \in
\mathcal{G}_{n}$, such that
%
\begin{equation}\label{DEFetatilde}
\| \widetilde{\eta}-\eta_{0}\| _{\infty}=O(
N_{n}^{-p}) .
\end{equation}
Let
%
\begin{equation} \label{DEFbetatilde}
\widetilde{\boldbeta}=\mathop{\arg\max}_{\tbbeta}
n^{-1}\sum_{i=1}^{n}Q[
g^{-1}\{ \widetilde{\eta}( \mathbf{X}_{i})
+\mathbf{Z}%
_{i}^{\T}\boldbeta\} ,Y_{i}] .
\end{equation}
In the following, let $m_{0i}\equiv m_{0}(
\mathbf{T}_{i}) =\eta_{0}( \mathbf{X}_{i})
+\mathbf{Z}_{i}^{\T}\boldbeta_{0}$ and $\varepsilon
_{i}=Y_{i}-g^{-1}(m_{0i})$. Further let
\[
\widetilde{m}_{0}( \mathbf{t}) =\widetilde{\eta}
( \mathbf{x}%
) +\mathbf{z}^{\T}\boldbeta_{0},\qquad \widetilde{m}%
_{0i}\equiv\widetilde{m}_{0}( \mathbf{T}_{i})
=\widetilde{\eta}%
( \mathbf{X}_{i}) +\mathbf{Z}_{i}^{\T}\boldbeta_{0}.
\]
\begin{lemma}\label{LEMbetatilde}
Under conditions \textup{(C1)--(C5)}, $\sqrt{n}(\widetilde{\boldbeta}-\boldbeta_{0})
\rightarrow\Normal(\mathbf{0},\break\mathbf{A}^{-1}\times\bolds\Sigma_{1}\mathbf{A}^{-1})$,
where {$\widetilde{\boldbeta}$ is in (\ref{DEFbetatilde})},
$\mathbf{A}=E[ \rho_{2}\{m_{0}( \mathbf{T}) \}\mathbf{Z}^{\otimes2}] $ and
$\bolds{\Sigma}_{1}=\break E[q_{1}^{2}\{ m_{0}{(\mathbf{T})} \} \mathbf{Z}^{\otimes2}]$.
\end{lemma}

In the following, denote $\widetilde{\boldtheta}= (
\widetilde{\boldgamma}^{\T}, \widetilde{\boldbeta}^{\T} )^{\T}$,
$\widehat{\boldtheta}=(\widehat{\boldgamma}^{\T},
\widehat{\boldbeta}^{\T})^{\T}$ and
%
\begin{equation} \label{DEFmtilde}
\widetilde{m}_{i}\equiv\widetilde{m}( \mathbf{T}_{i}) =
\widetilde{\eta}( \mathbf{X}_{i}) +\mathbf{Z}_{i}^{\T}
\widetilde{\boldbeta}=\mathbf{B}_{i}^{\T}\widetilde{\boldgamma}
+\mathbf{Z}_{i}^{\T}\widetilde{\boldbeta}.
\end{equation}
\begin{lemma}\label{LEMthetatilde-thetahat}
Under conditions \textup{(C1)--(C5)},
\[
\|
\widehat{\boldtheta}-\widetilde{\boldtheta}\|=O_{P}\{
N_{n}^{1/2-p}+(N_{n}/n)^{-1/2}\}.
\]
\end{lemma}

\subsection{\texorpdfstring{Proof of Theorem \protect\ref{THMconvergence}}{Proof of Theorem 1}}

According to Lemma \ref{LEMthetatilde-thetahat},
\begin{eqnarray*}
\| \widehat{\eta}-\widetilde{\eta}\| _{2}^{2}
&=&\| (
\widehat{\boldgamma}-\widetilde{\boldgamma}) ^{\T}\mathbf
{B}\| _{2}^{2}=(\widehat{%
\boldgamma}-\widetilde{\boldgamma}) ^{\T}E\Biggl[ n^{-1}%
\sum_{i=1}^{n}\mathbf{B}_{i}^{\otimes2} \Biggr] (
\widehat{\boldgamma}-
\widetilde{\boldgamma}) \\
&\leq&C\| \widehat{\boldgamma}-\widetilde{\boldgamma}%
\| _{2}^{2},
\end{eqnarray*}
thus $\| \widehat{\eta}-\widetilde{\eta}\| _{2}=
O_{P}\{ N_{n}^{1/2-p}+ (N_n/n)^{1/2}\}$ and
\begin{eqnarray*}
\| \widehat{\eta}-\eta_{0}\| _{2} &\leq&\|
\widehat{\eta}%
-\widetilde{\eta}\| _{2}+\| \widetilde{\eta}-\eta
_{0}\| _{2} = O_{P}\{ N_{n}^{1/2-p}+
(N_n/n)^{1/2}\}+O_{P}( N_{n}^{-p}) \\
&=& O_{P}\{ N_{n}^{1/2-p}+ (N_n/n)^{1/2}\}.
\end{eqnarray*}
By Lemma 1 of \citet{stone85}, $\| \widehat{\eta}_{k
}-\eta
_{0k }\| _{2k }= O_{P}\{ N_{n}^{1/2-p}+ (N_n/n)^{1/2}\}$,
for each $1\leq k \leq d_{1}$. Equation (\ref{DEFetatilde}) implies that
$\| \widehat{\eta}-\widetilde{\eta}\| _{n}= O_{P}\{
N_{n}^{1/2-p}+ (N_n/n)^{1/2}\}$. Then
\begin{eqnarray*}
\| \widehat{\eta}-\eta_{0}\| _{n} &\leq&\|
\widehat{\eta}%
-\widetilde{\eta}\| _{n}+\| \widetilde{\eta}-\eta
_{0}\| _{n}\\
&=& O_{P}\{N_{n}^{1/2-p}+(N_n/n)^{1/2}\}+O_{P}(
N_{n}^{-p}) \\
&=& O_{P}\{N_{n}^{1/2-p}+ (N_n/n)^{1/2}\}.
\end{eqnarray*}
Similarly,
\[
\sup _{\eta_{1},\eta_{2}\in\mathcal{S}_{n}^{0}}\biggl|
\frac{%
\langle\eta_{1},\eta_{2}\rangle_{n}-\langle
\eta_{1},\eta_{2}\rangle}{\| \eta_{1}\|
_{2}\| \eta_{2}\| _{2}}\biggr| =O_{P}\bigl\{ \bigl(
 \log(n) N_{n}/ n\bigr) ^{1/2}\bigr\}
\]
and $\| \widehat{\eta}_{k }-\eta_{0k }\| _{nk
}=O_{P}\{ N_{n}^{1/2-p}+ (N_n/n)^{1/2}\}$, for any $k =1,\ldots,d_1$.

\vspace*{15pt}\subsection{\texorpdfstring{Proof of Theorem \protect\ref{THMnormality}}{Proof of Theorem 2}}

We first verify that
%
\begin{eqnarray}
\label{EQIII1-1}
n^{-1}\sum_{i=1}^{n}\rho_{2}( m_{0i})
\widetilde{\mathbf{Z}}_{i}\Gamma(\mathbf{X}_{i})
^{\T}(
\widehat{\boldbeta}-\boldbeta_{0}) &=&o_{P}(
n^{-1/2}), \\
\label{EQIII1-2}
n^{-1}\sum_{i=1}^{n}\{(\widehat{\eta}-\eta_{0})(
\mathbf{X}_{i}) \} \rho_{2}( m_{0i})
\widetilde{\mathbf{Z}}_{i}&=&o_{P}( n^{-1/2}),
\end{eqnarray}
where $\widetilde{\mathbf{Z}}$ is defined in (\ref{DEFGammax}).

Define
%
\begin{equation} \label{DEFMn}
\mathcal{M}_{n}=\{ m( \mathbf{x},\mathbf{z})
=\eta
( \mathbf{x}) +\mathbf{z}^{\T}\bolds\boldbeta\dvtx\eta\in
\mathcal{G}_{n}\}.
\end{equation}
Noting that $\rho_{2}$ is a fixed bounded function under (C7), we
have $E[(\widehat{\eta}-\eta_{0})( \mathbf{X}) \rho
_{2}( m_{0})
\widetilde{Z}_{l}]^2 \leq O( \| \widehat{m}%
-m_{0}\| _{2}^{2})$, for $l=1,\ldots,d_{2}$. By Lemma
\ref{LEMEntropy}, the logarithm of the $\varepsilon$-bracketing
number of the class of functions
\[
\mathcal{A}_{1}( \delta)
=\bigl\{\rho_{2}\{m(\mathbf{x},\mathbf{z})\}
\{\mathbf{z}-\Gamma(\mathbf{x})\}\dvtx m\in\mathcal{M}_{n},\|
m-m_{0}\| \leq\delta\bigr\}
\]
is $c\{ N_{n} \log(\delta/\varepsilon) +\log
( \delta^{-1}) \} $, so the corresponding
entropy integral
\[
\mathcal{J}_{[ \cdot] }( \delta
,\mathcal{A}_{1}( \delta) ,\| \cdot\| )
\leq c\delta\{ N_{n} ^{1/2}+ \log^{1/2}( \delta
^{-1}) \} .
\]
According to Lemmas \ref{LEMbetatilde} and
\ref{LEMthetatilde-thetahat} and Theorem \ref{THMconvergence},
$\|\widehat{m}-m_{0}\|_{2}=O_{P}\{
N_{n}^{1/2-p}+(N_{n}/n)^{1/2}\}$.
By Lemma 7 of \citet{stone86}, we have $ \| \widehat{\eta
}-\eta
_{0}\| _{\infty}\leq cN_{n}^{1/2}\| \widehat{\eta
}-\eta_{0}\| _{2}=O_{P}(N_{n}^{1-p}+N_{n} n^{-1/2})$, thus
%
\begin{equation}\label{EQmhat-m-infnorm}
\| \widehat{m}-m_{0}\| _{\infty}=
O_{P}(N_{n}^{1-p}+N_{n} n^{-1/2} ).
\end{equation}
Thus by Lemma \ref{LEMExpection-entropy} and Theorem
\ref{THMconvergence}, for $r_{n}=\{N_{n}^{1/2-p}+
(N_{n}/n)^{1/2}\}^{-1}$,
\begin{eqnarray*}
&&E\Biggl| n^{-1}\sum_{i=1}^{n}\{(\widehat{\eta}-\eta
_{0})( \mathbf{X}_{i}) \} \rho_{2}(
m_{0i}) \widetilde{\mathbf{Z}}_{i}-E[ (\widehat{\eta
}-\eta_{0})( \mathbf{X}) \rho_{2}\{ m_{0}(
\mathbf{T}) \}
\widetilde{\mathbf{Z}} ] \Biggr| \\
&&\qquad\leq n^{-1/2}Cr_{n}^{-1} \{ N_{n} ^{1/2}+
\log^{1/2}( r_{n} ) \} \biggl[ 1+\frac{cr_{n}^{-1}
\{N_{n}^{1/2}+\log^{1/2}(r_{n})\}
}{r_{n}^{-2}\sqrt{n}}
M_{0}\biggr] \\
&&\qquad\leq
O(1)n^{-1/2}r_{n}^{-1}\{N_{n}^{1/2}+\log^{1/2}(r_{n})\},
\end{eqnarray*}
where\vspace*{2pt}
$r_{n}^{-1}\{N_{n}^{1/2}+\log^{1/2}(r_{n})\}=o(
1)$ according to condition (C5). By the definition of
$\widetilde{\mathbf{Z}}$, for any measurable function $\phi$,
$E[ \phi( \mathbf{X})\rho_{2}\{ m_{0}
( \mathbf{T}%
) \} \widetilde{\mathbf{Z}}] =\mathbf{0}$. Hence
(\ref{EQIII1-2}) holds. Similarly, (\ref{EQIII1-1}) follows from
Lemmas \ref{LEMExpection-entropy} and
\ref{LEMthetatilde-thetahat}.

According to condition (C6), the projection function $\Gamma
^{\mathrm{add}}( \mathbf{x}) =\break\sum_{k=1}^{d_{1}}\Gamma
_{k}(x_{k})$, where the theoretically centered function $\Gamma
_{k}\in\mathcal{H}{(p)}$. By the result
of de Boor [(\citeyear{deBoor01}), page 149], there exists an
empirically centered function $%
\widehat{\Gamma}_{k}\in\mathcal{S}_{n}^{0}$, such that $\|
\widehat{\Gamma}_{k}-{\Gamma}_{k}\| _{\infty
}=O_{P}(N_{n}^{-p})$, $k=1,\ldots,d_{1}$. Denote $\widehat{\Gamma
}^{\mathrm{add}}( \mathbf{x})
=\sum_{k=1}^{d_{1}}\widehat{\Gamma}_{k}(x_{k})$ and clearly
$\widehat{%
\Gamma}^{\mathrm{add}}\in\mathcal{G}_{n}$. For any $\boldnu\in
R^{d_{2}}$,
define $\widehat{m}_{\tbnu}=\widehat{m}( \mathbf{x},\mathbf
{z}) +%
\boldnu^{\T}\{ \mathbf{z}-\widehat{\Gamma}^{\mathrm{add}}(
\mathbf{x}%
) \} =\{ \widehat{\eta}( \mathbf{x}
) -\boldnu^{%
\T}\widehat{\Gamma}^{\mathrm{add}}( \mathbf{x}) \}
+( \widehat{\boldbeta}+\boldnu) ^{\T}\mathbf{z}\in
\mathcal{M}_{n}$,
where $\mathcal{M}_{n}$ is given in (\ref{DEFMn}). Note that
$\widehat{m}%
_{\tbnu}$ maximizes the function $\widehat{l}_{n}( m)
=n^{-1}\sum_{i=1}^{n}Q[ g^{-1}\{ m(
\mathbf{T}_{i})
\} ,Y_{i}] $ for all $m\in\mathcal{M}_{n}$ when
$\boldnu=%
\mathbf{0}$, thus\break $ \frac{\partial}{\partial\tbnu}\widehat
{l}_{n}( \widehat{m}_{\tbnu}) |_{\tbnu=\mathbf{0}}=\mathbf{0}$. For simplicity, denote
$\widehat{m}_{i}\equiv\widehat{m}( \mathbf{T}_{i})$,
and we have
%
\begin{equation}\label{EQ1stdiv}
\mathbf{0}\equiv\frac{\partial}{\partial\tbnu}\widehat{l}%
_{n}( \widehat{m}_{\tbnu}) \bigg| _{\tbnu=\mathbf{0}}
=n^{-1}\sum_{i=1}^{n}q_{1}( \widehat{m}_{i},Y_{i})
\widetilde{\mathbf{Z}}_{i} +O_{P}(N_{n}^{-p}).
\end{equation}
For the first term in (\ref{EQ1stdiv}), we get
%
\begin{eqnarray}\label{EQdecomp}
n^{-1}\sum_{i=1}^{n}q_{1}( \widehat{m}_{i},Y_{i})
\widetilde{\mathbf{Z}}_{i}&=&n^{-1}\sum_{i=1}^{n}q_{1}(
m_{0i},Y_{i}) \widetilde{\mathbf{Z}}_{i}\nonumber\\
&&{}+n^{-1}\sum_{i=1}^{n}q_{2}( m_{0i},Y_{i})
( \widehat{m}%
_{i}-m_{0i}) \widetilde{\mathbf{Z}}_{i}\nonumber\\[-8pt]\\[-8pt]
&&{} +n^{-1}\sum_{i=1}^{n}q_{2}^{\prime}(
\bar{m}_{i},Y_{i}) ( \widehat{m}_{i}-m_{0i})
^{2}\widetilde{\mathbf{Z}}_{i}\nonumber\\
&=&\mathrm{I}+\mathrm{II}+\mathrm{III}.\nonumber
\end{eqnarray}
We decompose $\mathrm{II}$ into two terms $\mathrm{II}_{1}$ and
$\mathrm{II}_{2}$ as follows:
\begin{eqnarray*}
\mathrm{II}&=&n^{-1}\sum_{i=1}^{n}q_{2}( m_{0i},Y_{i})
\widetilde{\mathbf{Z}}_{i} \{ (\widehat{\eta}-\eta_{0})
( \mathbf{X}%
_{i}) \}+n^{-1}\sum_{i=1}^{n}q_{2}(
m_{0i},Y_{i})
\widetilde{\mathbf{Z}}_{i} \mathbf{Z}_{i}^{\T}( \widehat
{\boldbeta}-\boldbeta%
_{0}) \\
&=&\mathrm{II}_{1}+\mathrm{II}_{2}.
\end{eqnarray*}
We next show that
%
\begin{equation}\label{EQII1}
\mathrm{II}_{1}=\mathrm{II}_{1}^{*}+o_{P}(
n^{-1/2}),
\end{equation}
where $\mathrm{II}_{1}^{*}=-n^{-1}\sum_{i=1}^{n}\rho_{2}(
m_{0i})
\widetilde{\mathbf{Z}}_{i} \{ (\widehat{\eta}-\eta_{0})
( \mathbf{X}_{i}) \}$. Using an argument similar to the proof
of Lemma \ref{LEMthetatilde-thetahat}, we have
\[
(\widehat{\eta}-\eta_{0})( \mathbf{X}_{i}) =\mathbf
{B}_{i}^{\T%
}{\mathbf{K}}\mathbf{V}_{n}^{-1}\Biggl\{
n^{-1}\sum_{i=1}^{n}q_{1}( m_{0i},Y_{i})
\mathbf{D}_{i}^{\T}+o_{P}( N_{n}^{-p}) \Biggr\},
\]
where\vspace*{1pt} $\mathbf{K}=(\mathbf{I}_{N_nd_1},\mathbf{0}_{(N_nd_1)\times
d2})$ {and} $\mathbf{I}_{N_nd_1}$ is a diagonal matrix. Note that the
expectation of the square of the $s$th column of
$n^{-1/2}(\mathrm{II}_{1}-\mathrm{II}_{1}^{*})$ is
\begin{eqnarray*}
&&E\Biggl[ n^{-1/2}\sum_{i=1}^{n}\{ q_{2}(
m_{0i},Y_{i}) +\rho_{2}( m_{0i}) \}
\widetilde{Z}_{is}(\widehat{\eta}-\eta_{0})(
\mathbf{X}_{i}) \Biggr] ^{2} \\
&&\qquad= n^{-1}\sum_{i=1}^{n}\sum_{j=1}^{n}E\{
\varepsilon_{i}\varepsilon_{j}\rho_{1}^{\prime}(
m_{0i}) \rho_{1}^{\prime}(
m_{0{j}})
\widetilde{Z}_{is}\widetilde{Z}_{js}(\widehat{\eta}-\eta
_{0})( \mathbf{X}_{i}) (\widehat{\eta}-\eta
_{0})(
\mathbf{X}_{j}) \}\\
&&\qquad=n^{-3}\sum_{i=1}^{n}\sum_{j=1}^{n}\sum_{k=1}^{n}\sum
_{l=1}^{n}E\{
\varepsilon_{i}\varepsilon_{j}\varepsilon_{k}\varepsilon
_{l}\rho_{1}^{\prime}( m_{0i}) \rho_{1}^{\prime
}( m_{0j})
\rho_{1}( m_{0k}) \rho_{1}( m_{0l}) \\
&&\qquad\quad\hphantom{=n^{-3}\sum_{i=1}^{n}\sum_{j=1}^{n}\sum_{k=1}^{n}\sum
_{l=1}^{n}}\hspace*{29.4pt}{}\times\widetilde{Z}_{is}\widetilde{Z}_{js}\mathbf
{B}_{i}^{\T}\mathbf{K}\mathbf{V}%
_{n}^{-1}\mathbf{D}_{i}^{\T}\mathbf{B}_{j}^{\T}\mathbf{K}\mathbf
{V}%
_{n}^{-1}\mathbf{D}_{j}^{\T}\}\\
&&\qquad\quad{} + o(nN_n^{-2p}) = o(
1), \qquad s=1,\ldots,d_2.
\end{eqnarray*}
Thus, (\ref{EQII1}) holds by Markov's inequality.
Based on (\ref{EQIII1-2}), we have $ \mathrm{II}_{1}^{*}
=o_{P}( n^{-1/2})$. Using similar arguments and
(\ref{EQIII1-1}) and (\ref{EQIII1-2}), we can show that
\begin{eqnarray*}
\mathrm{II}_{2} &=& -n^{-1}\sum_{i=1}^{n}\rho_{2}(
m_{0i}) \widetilde{\mathbf{Z}}_{i} \mathbf{Z}_{i}^{\T}(
\widehat{\boldbeta}
-\boldbeta_{0}) +o_{P}( n^{-1/2}) \\
&=&-n^{-1}\sum_{i=1}^{n}\rho_{2}( m_{0i})
\widetilde{\mathbf{Z}}_{i}^{\otimes2}( \widehat{\boldbeta}-
\boldbeta_{0}) +o_{P}( n^{-1/2}).
\end{eqnarray*}
According to (\ref{EQmhat-m-infnorm}) and condition (C5), we have
\begin{eqnarray*}
\mathrm{III} &=&n^{-1}\sum_{i=1}^{n}q_{2}^{\prime}(
\bar{m}_{i},Y_{i})
( \widehat{m}_{i}-m_{0i}) ^{2}\widetilde{\mathbf
{Z}}_{i} \\
&\leq&C\| \widehat{m}-m_{0}\| _{\infty
}^{2}=O_{p}\bigl\{ N_{n}^{2( 1-p)}
+N_{n}^{2}n^{-1}\bigr\}\\
&=&o_{P}( n^{-1/2}).
\end{eqnarray*}
Combining (\ref{EQ1stdiv}) and (\ref{EQdecomp}), we have
\[
\mathbf{0}=n^{-1}\sum_{i=1}^{n}q_{1}( m_{0i},Y_{i})
\widetilde{\mathbf{Z}}_{i} +\bigl\{ E[ \rho_{2}\{
m_{0}( \mathbf{T}) \}
\widetilde{\mathbf{Z}}^{\otimes2}]+o_{P}( 1)
\bigr\} (\widehat{\boldbeta}-\boldbeta_{0})+o_{P}(
n^{-1/2}).
\]
Note that
\[
E[ \rho_{1}^{2}\{m_{0} (\mathbf{T})\}
\varepsilon^{2}
\widetilde{\mathbf{Z}}^{\otimes2}]=E[E(\varepsilon
^{2}|\mathbf{T})\rho_{1}^{2}\{m_{0}
(\mathbf{T})\}\widetilde{\mathbf{Z}}^{\otimes
2}]=E[
\rho_{2}\{ m_{0}( \mathbf{T}) \}
\widetilde{\mathbf{Z}}^{\otimes2}].
\]
Thus the desired distribution of $\widehat{\boldbeta}$ follows.

\subsection{\texorpdfstring{Proof of Theorem \protect\ref{THMrootn}}{Proof of Theorem 3}}
Let $\tau_{n}=n^{-1/2}+a_{n}$. It suffices to show that for any
given $\zeta>0$, there exists a large constant $C$ such that
%
\begin{equation}\label{f2}
\pr\Bigl\{\sup_{\Vert\mathbf{u}\Vert
=C}\mathcal{L}_{P}(\boldbeta_{0}+\tau
_{n}\mathbf{u})<\mathcal{L}_{P}(\boldbeta_{0})\Bigr\} \geq1-\zeta.
\end{equation}
Denote
\begin{eqnarray*}
U_{n,1} &=&\sum_{i=1}^{n}\bigl[ Q\bigl\{ g^{-1}\bigl(
\widehat{\eta}^{\mpl} (\mathbf{X}_{i})+\mathbf{Z}_{i}^{\T}(
\boldbeta_{0}+\tau_{n}\mathbf{u}) \bigr) ,Y_{i}\bigr\}
\\
&&\hspace*{42.6pt}{}-Q\bigl\{
g^{-1}\bigl(\widehat{\eta}^{\mpl}(\mathbf{X}_{i})+\mathbf
{Z}_{i}^{\T}\boldbeta_{0}\bigr)
,Y_{i}\bigr\} \bigr]
\end{eqnarray*}
and $U_{n,2}=-n\sum_{j=1}^{s}\{p_{\lambda_{n}}(|\beta
_{j0}+\tau_{n}v_{j}|)-p_{\lambda_{n}}( | \beta _{j0}| ) \}$, where $s$
is the number of components of $\bolds{\beta}_{10}$. Note that $
p_{\lambda_{n}}( 0) =0$ and $p_{\lambda_{n}}( | \beta| ) \ge0$ for all
$\beta$.\vspace*{1pt} Thus, $
\mathcal{L}_{P}(\boldbeta_{0}+\tau_{n}\mathbf{u})-\mathcal{L}_{P}(%
\boldbeta_{0})\le U_{n,1}+U_{n,2}$. Let
$\widehat{m}_{0i}^{\mpl}=\widehat{\eta}^{\mpl}(\mathbf
{X}_{i})+\mathbf{Z}_{i}^{\T}\boldbeta_{0}$.
For~$U_{n,1}$, note that
\[
U_{n,1}=\sum_{i=1}^{n}[ Q\{ g^{-1}(
\widehat{m}_{0i}^{\mpl}+\tau_{n}\mathbf{u}^{\T}\mathbf
{Z}_{i})
,Y_{i}\} -Q\{ g^{-1}(
\widehat{m}_{0i}^{\mpl}) ,Y_{i}\} ] .
\]
Mimicking the proof for Theorem \ref{THMnormality} indicates that
%
\begin{equation}\label{EQDn1}
U_{n,1}=\tau
_{n}\mathbf{u}^{\T}\sum_{i=1}^{n}q_{1}(m_{0i},Y_{i})\widetilde
{\mathbf{Z}}
_{i}+{\frac{n}{2}}\tau
_{n}^{2}\mathbf{u}^{\T}\bolds{\Omega}\mathbf{u}+o_{P}(1),
\end{equation}
where the orders of the first term and the second term are $%
O_{P}(n^{1/2}\tau_{n}) $ and $O_{P}(n\tau_{n}^{2})$, respectively.
For $%
U_{n,2}$, by a Taylor expansion and the Cauchy--Schwarz inequality, $%
n^{-1}U_{n,2}$ is bounded by $ \sqrt{s}\tau_{n}a_{n}\Vert
\mathbf{u}\Vert+\tau_{n}^{2}w_{n}\Vert\mathbf{u}\Vert
^{2}=\break C\tau_{n}^{2}(\sqrt{s}+w_{n}C)$.
As $w_{n}\to0$, both the first and second terms on the right-hand side
of (%
\ref{EQDn1}) dominate $U_{n,2}$, by taking $C$ sufficiently large.
Hence, (%
\ref{f2}) holds for sufficiently large $C$.

\subsection{\texorpdfstring{Proof of Theorem \protect\ref{THMoracle}}{Proof of Theorem 4}}

The proof of $\widehat{\boldbeta}_{2}^{\mpl}=0$ is similar to that
of Lem\-ma~3 in \citet{liliang08}. We therefore omit the details
and refer to the proof of that lemma.

Let $\widehat{m}^{\mpl}(\mathbf{x},\mathbf{z}_1)=\widehat{\eta
}^{\mpl}(%
\mathbf{x})+\mathbf{z}_{1}^{\T}\boldbeta_{10}$, for
$\widehat{\eta}^{\mpl}$ in
(\ref{DEFetaMPL}), and $m_{0}(\mathbf{T}_{1i})=%
\bolds{\eta}_{0}^{\T}(\mathbf{X}_{i})+\mathbf{Z}_{i1}^{\T
}\boldbeta%
_{10}$.
Define $\mathcal{M}_{1n}=\{ m(
\mathbf{x},\mathbf{z}_{1}) =\eta
( \mathbf{x}) +\mathbf{z}_{1}^{\T}\boldbeta_{1}\dvtx\eta
\in\mathcal{G}%
_{n}\}$.
For any $\boldnu_{1}\in R^{s}$, where $s$ is the dimension of
$\boldbeta_{10}$%
, define
\[
\widehat{m}_{\tbnu_{1}}^{\mpl}(\mathbf{t}_{1})=\widehat{m}(
\mathbf{x},\mathbf{z}%
_1) +\boldnu_{1}^{\T}\widetilde{\mathbf{z}}_{1}=\{
\widehat{\eta}^{\mpl}( \mathbf{x}) -%
\boldnu_{1}^{\T}\Gamma_{1}( \mathbf{x}) \}
+(\widehat{\boldbeta}_{1}^{\mpl}+\boldnu_{1}) ^{\T}\mathbf{z}_{1}.
\]
Note that\vspace*{1pt} $\widehat{m}_{\tbnu_{1}}^{\mpl}$ maximizes
$\sum_{i=1}^{n}Q[ g^{-1}\{ m_{0}(\mathbf{T}_{1i})\},Y_{i}]-
n\sum_{j=1}^{s}p_{\lambda_{jn}}(|\widehat{\beta} _{j1}^{\mpl}+\break v_{j1}|)$
for all $m \in\mathcal{M}_{1n}$ when $\boldnu_{1}=\mathbf{0}$. Mimicking the proof for Theorem
\ref{THMnormality} indicates that
\begin{eqnarray*}
\mathbf{0}&=&n^{-1}\sum_{i=1}^{n}q_{1}\{
m_{0}(\mathbf{T}_{1i}),Y_{i}\}
\widetilde{\mathbf{Z}}_{1i}+\{ p_{\lambda_{jn}}^{\prime
}( | \beta_{j0}| ) \operatorname{sign}
( \beta
_{j0}) \}
_{j=1}^{s}+o_{P}( n^{-1/2})  \\
&&{}+\bigl\{ E[ \rho_{2}\{ m_{0}(
\mathbf{T}_{1}) \} \widetilde{\mathbf{Z}}_{1}^{\otimes
2}] +o_{P}( 1) \bigr\} (
\widehat{\boldbeta}_{1}^{\mpl}-\boldbeta_{10}) \\
&&{}+\Biggl\{ \sum_{j=1}^{s}p_{\lambda
_{jn}}^{\prime\prime}(
| \beta_{j0}| ) +o_{P}( 1) \Biggr\}
(\widehat{\bolds{\beta}}_{j1}^{\mpl}-\beta_{j0}).
\end{eqnarray*}
Thus, asymptotic normality follows because
\begin{eqnarray*}
\mathbf{0}&=&n^{-1}\sum_{i=1}^{n}q_{1}\{
m_{0}(\mathbf{T}_{1i}),Y_{i}\}
\widetilde{\mathbf{Z}}_{1i}+\xi_{n}+o_{P}(
n^{-1/2})\\
&&{}+\{ \bolds{\Omega}_{s}+\bolds{\Sigma
}_{\lambda}+o_{P}( 1) \} (
\widehat{\boldbeta}_{1}^{\mpl}-\boldbeta_{10}),\\[-22pt]
\end{eqnarray*}
\[
E[ \rho_{1}^{2}\{m_{0}
(\mathbf{T}_{1})\} \{Y-m_{0}
(\mathbf{T}_{1})\}^{2}
\widetilde{\mathbf{Z}}_{1}^{\otimes2}]=E[ \rho
_{2}\{ m_{0}( \mathbf{T}_{1}) \}
\widetilde{\mathbf{Z}}_{1}^{\otimes2}].
\]
\end{appendix}

\section*{Acknowledgments}

The authors would like to thank the Co-Editor, Professor Tony Cai, an
Associate Editor and two referees for their constructive comments that
greatly improved an earlier version of this paper.

\begin{supplement}[id=suppA]
\stitle{Detailed proofs and additional simulation results of:
Estimation and variable selection for generalized additive partial
linear models\\}
\slink[doi]{10.1214/11-AOS885SUPP} 
\sdatatype{.pdf}
\sfilename{aos885\_suppl.pdf}
\sdescription{The supplemental materials contain detailed proofs and additional simulation results.}
\end{supplement}


%
\printaddresses

\end{document}